\documentclass[3p,preprint]{elsarticle}

\usepackage[utf8]{inputenc}
\usepackage[english]{babel}
\usepackage{amsfonts}
\usepackage{amsmath}
\usepackage{amssymb}
\usepackage{amsbsy}
\usepackage{graphicx}
\usepackage{float}
\usepackage{xcolor}
\usepackage{bbm}
\usepackage{gensymb}
\usepackage{subfigure}
\usepackage{multirow}
\usepackage{enumitem}
\usepackage{caption}
\usepackage{adjustbox}
\usepackage{comment}
\usepackage{algorithm}
\usepackage[noend]{algpseudocode}
\usepackage{url}

\newcommand{\xx}{\boldsymbol}

\newcommand{\ww}{\widehat}
\makeatletter
\newcommand{\xleftrightarrow}[2][]{%
	\ext@arrow3399{\longleftrightarrowfill@}{#1}{#2}}
\newcommand{\longleftrightarrowfill@}{%
	\arrowfill@\leftarrow\relbar\rightarrow}
\makeatother

\journal{Computer Methods in Applied Mechanics and Engineering }

\begin{document}
\setcounter{page}{1}
\begin{frontmatter}



\title{Modeling cardiac muscle fibers in ventricular and atrial electrophysiology simulations} 


\author[a]{Roberto Piersanti}
\ead{roberto.piersanti@polimi.it}

\author[a]{Pasquale C. Africa}
\ead{pasqualeclaudio.africa@polimi.it}

\author[a]{Marco Fedele}
\ead{marco.fedele@polimi.it}

\author[b]{Christian Vergara\corref{cor}}
\ead{christian.vergara@polimi.it}
\cortext[cor]{Corresponding author}

\author[a]{Luca Ded{\`e}}
\ead{luca.dede@polimi.it}

\author[c]{Antonio F. Corno}
\ead{antonio.f.corno@uth.tmc.edu}

\author[a,d]{Alfio Quarteroni} 
\ead{alfio.quarteroni@polimi.it}

\address[a]{Modellistica e Calcolo Scientifico (MOX), Dipartimento di Matematica, Politecnico di Milano, Milan, Italy}
\address[b]{Laboratory of Biological Structure Mechanics (LaBS), Dipartimento di Chimica, Materiali e Ingegneria Chimica, Politecnico di Milano, Milan, Italy}

\address[c]{\textcolor{black}{Houston Children's Heart Institute, Hermann Children's Hospital, University of Texas Health, McGovern Medical School, Houston, Texas, USA}}

\address[d]{Institute of Mathematics, \'Ecole Polytechnique Fédérale de Lausanne, Switzerland (Professor emeritus)}

\begin{abstract}
\noindent
Since myocardial fibers drive the electric signal propagation throughout the myocardium, accurately modeling their arrangement is essential for simulating heart electrophysiology (EP). Rule-Based-Methods (RBMs) represent a commonly used strategy to include cardiac fibers in computational models. A particular class of such methods is known as Laplace-Dirichlet-Rule-Based-Methods (LDRBMs) since they rely on the solution of Laplace problems. In this work we provide a unified framework, based on LDRBMs, for generating full heart muscle fibers. First, we \textcolor{black}{review existing ventricular LDRBMs providing a communal mathematical description} 
and introducing also some modeling improvements with respect to the existing literature. We then carry out a systematic comparison of LDRBMs based on meaningful biomarkers produced by numerical EP simulations. Next we propose, for the first time, a LDRBM to be used for generating atrial fibers. The new method, tested both on idealized and realistic atrial models, can be applied to any arbitrary geometries. Finally, we present numerical results obtained in a realistic whole heart where fibers are included for all the four chambers using the discussed LDRBMs.
\end{abstract}




\end{frontmatter}


\section{Introduction}\label{sec:intro}
In numerical heart electrophysiology (EP) a critical issue is that of modeling the myocardial fibers arrangement that characterizes the cardiac tissue. Aggregations of myofibers, namely the results of cardiomyocytes orientation, determine how the electric signal propagates within the muscle \cite{streeter1969fiber,roberts1979influence,punske2005effect}. \textcolor{black}{Moreover, also the muscle mechanical contraction, which is triggered by the propagation of the electric signal thorough the tissue, strongly depends on the fibers orientation \cite{papadacci2017imaging,gil2019influence,eriksson2013influence,palit2015computational,guan2020effect}}. This motivates the need to accurately include fiber orientations in order to obtain physically meaningful results \cite{beyar1984computer,bayer2005laplace}.

In the last decades, myofibers orientation have been studied using histological data and Diffusion Tensor Imaging (DTI) acquisitions \cite{hsu1998magnetic,helm2005ex,pashakhanloo2016myofiber}. DTI is a Magnetic Resonance Imaging (MRI) technique able to produce useful structural information about heart muscle fibers and
largely applied to explanted ex-vivo hearts, coming from animal experiments \cite{hsu1998magnetic,helm2005ex,scollan2000reconstruction,wu2007study,peyrat2007computational} or from human corpses~ \cite{pashakhanloo2016myofiber,lombaert2012human}. However, acquired in-vivo DTI protocol lasts hours and generally produces a noisy low-resolution fibers reconstruction~\cite{toussaint2013vivo,alexander2001analysis,nagler2013personalization}. Furthermore, since the \textcolor{black}{atrial} thickness is smaller than the DTI voxel size, it is not possible to obtain in-vivo myofibers in the atria \cite{hoermann2019automatic}. All the above considerations make nowadays DTI technique unusable to reconstruct accurate 3D myofibers field in the common clinical practice. 

Because of the difficulties to acquire patient-specific fibers data, different methodologies have been proposed in order to provide a realistic surrogate of fiber orientation for in-vivo cardiac geometries \cite{lombaert2012human,hoermann2019automatic,bayer2012novel,wong2014generating,rossi2014thermodynamically,doste2019rule,krueger2011modeling,tobon2013three,fastl2018personalized}. Among these, atlas-based methods map and project a detailed fiber field, previously reconstructed on an atlas, on the geometry of interest , exploiting DTI or histological data; see \cite{lombaert2012human} for the ventricles and \cite{hoermann2019automatic,roney2020constructing} for the atria. However, these methods require complex registration algorithms and are strictly dependent on the original atlas data upon which they have been built. 

Alternative strategies for generating myofiber orientations are the so called {\sl Rule-Based Methods} (RBMs) \cite{beyar1984computer,potse2006comparison,nielsen1991mathematical,bishop2009development}. RBMs describe fiber orientations with mathematically sound rules based on histological or DTI observations and require information only about the myocardial geometry. These methods parametrize the transmural and apico-basal directions in the entire myocardium in order to assign orthotropic (longitudinal, transversal and normal) myofibers; see \cite{bayer2005laplace,bayer2012novel,wong2014generating,rossi2014thermodynamically,doste2019rule,potse2006comparison} for the ventricles and \cite{krueger2011modeling,tobon2013three,fastl2018personalized,ferrer2015detailed,plank2008evaluating} for the atria.

A particular class of RBMs, which rely on the solution of Laplace boundary-value problems, is known as {\sl Laplace-Dirichlet-Rule-Based-Methods} (LDRBMs), addressed in \cite{bayer2005laplace,bayer2012novel,wong2014generating,rossi2014thermodynamically,doste2019rule} for the ventricular case. LDRBMs define the transmural and apico-basal directions by taking the gradient of solutions corresponding to suitable Dirichlet boundary conditions. These directions are then properly rotated in order to match histological observations \cite{greenbaum1981left,sanchez1990myocardial,anderson2009three}. The above procedure ensures a smooth and continuous change in fibers directions throughout the whole myocardium. 

Most of existing ventricular RBMs refer to left ventricle only and usually introduce an artificial basal plane located well below the cardiac valves. Only recently, a LDRBM, that takes into account fiber directions in specific cardiac regions, such as the right ventricle, the inter-ventricular septum and the outflow tracks, has been developed \cite{doste2019rule}. This has provided a great improvement in RBMs since the right ventricle exhibits a different fiber orientation with respect to the left ventricle \cite{helm2005ex,scollan2000reconstruction,lombaert2012human,kocica2006helical}. The presence of a discontinuity in the inter-ventricular septal fibers is a crucial matter, still very debated \cite{kocica2006helical,boettler2005new}, even though the corresponding effects on electrical signal propagation have not been studied yet.

Regarding the atria, several RBMs have been developed. They either use semi-automatic approaches \cite{krueger2011modeling,tobon2013three,fastl2018personalized,ferrer2015detailed,plank2008evaluating,krueger2010patient,labarthe2012semi,rocher2019highly} or prescribe manually the fiber orientations in specific atrial regions \cite{harrild2000computer,vigmond2001reentry,jacquemet2003study,seemann2006heterogeneous}. Often, these procedures require a manual intervention and, in many cases, are designed for specific atrial morphologies \cite{krueger2011modeling,fastl2018personalized}. Hence, a general automatic processing pipeline for generating atrial fibers field still remains a knotty
procedure \cite{fastl2018personalized,fastl2016personalized}. Moreover, no LDRBMs have been proposed so far for the atria. As a matter of fact, an extension of the ventricular LDRBMs is not straightforward, mainly because the atrial fibers architecture is characterized by the presence of multiple overlapping bundles running along different directions, differently from the ventricles
one where myofibers are aligned along regular patterns.

Over the past years several cardiac computational models were carried out in order to study pathological conditions affecting either the electrical or mechanical response in individual heart chambers. However, in the quest for a more quantitative understanding of the heart functioning both in health and diseased scenarios, it became fundamental to model and simulate the entire heart as an whole organ  \cite{trayanova2011whole}. Only recently, the scientific community moved towards the whole heart modeling and simulations \cite{baillargeon2014living,fritz2014simulation, augustin2016anatomically,quarteroni2017integrated,santiago2018fully,land2018influence,pfaller2019importance,strocchi2020simulating}. Nevertheless, we highlight that
none of these computational studies makes use of a unified methodology to embed reliable and detailed cardiac fibers in the whole heart muscle to take into account different fiber orientations specific of the four chambers.  

Driven by the previous open issues, in this work we provide a unified framework, based on LDRBMs, for generating full heart muscle fibers. We start by giving a 
\textcolor{black}{review based on a communal mathematical} description for three existing LDRBMs in the ventricles, introducing also some modeling improvements \cite{bayer2012novel,rossi2014thermodynamically,doste2019rule}. In particular, we extend ventricular LDRBMs in order to include specific fiber directions for the right ventricle~\textcolor{black}{\footnote{\textcolor{black}{Notice that the original LDRBM by Doste et al. \cite{doste2019rule} already includes specific fiber directions for the right ventricle.}}}. Next, we carry out a systematic comparison of the effect produced by different LDRBMs on EP in terms of meaningful biomarkers (e.g. activation times) computed from numerical simulations. Specifically, we study the importance of including different fiber orientations in the right ventricle and we investigate the effect of the inter-ventricular septal fibers discontinuity. 

Then, at the best of our knowledge, we propose for the first time an atrial LDRBM which is able to qualitatively reproduce all the important features, such as fiber bundles, needed to provide a realistic atrial musculature architecture. Unlike most of the existing RBMs, the new method, tested both on idealized and realistic atrial geometries, can be easily applied to any arbitrary geometries. 
Moreover, we analyse the strong effect of the complex atrial fiber architecture on the electric signal propagation obtained by numerical simulations.

In the last part of the work, we illustrate numerical results including the full heart LDRBMs fiber generations and an EP simulation with physiological activation sites in a four chamber realistic computational domain of the heart.

This paper is organized as follows. In Section \ref{sec:ventr} we review and provide a 
\textcolor{black}{communal mathematical} description of existing LDRBMs for ventricular fiber generation. In Section \ref{sec:atrial} we propose the novel LDRBM for atrial fibers generation. Numerical methods to perform EP simulations are explained in Section \ref{sec:num}. Section~\ref{sec:res} is dedicated to numerical results where we present a comparison among different ventricular LDRBMs, we test the new atrial LDRBM and we show numerical simulation of the whole heart EP including the presented fiber generation methods. \textcolor{black}{Section \ref{sec:disc} is devoted to the discussion of the results; in particular, we compare the results of the novel atrial LDRBM with the fiber orientations obtained by the RBM proposed in \cite{tobon2013three,ferrer2015detailed} and also to anatomical data \cite{sanchez2014left}. Moreover, we discuss the limitations of the presented study.} Finally, conclusions follow.  

\section{Rule-Based-Methods for ventricular fibers generation}\label{sec:ventr}

In this section we review three popular LDRBMs introduced so far in the literature: specifically, we consider LDRBMs by Rossi et al. \cite{rossi2014thermodynamically}, by Bayer et al. \cite{bayer2012novel} and by Doste et al. \cite{doste2019rule}. In view of our review process of the former LDRBMs, we provide a 
\textcolor{black}{communal mathematical} framework of such methods, highlighting similarities and differences. 

\begin{figure}[t!]
	\centering
	\includegraphics[width=0.85\textwidth]{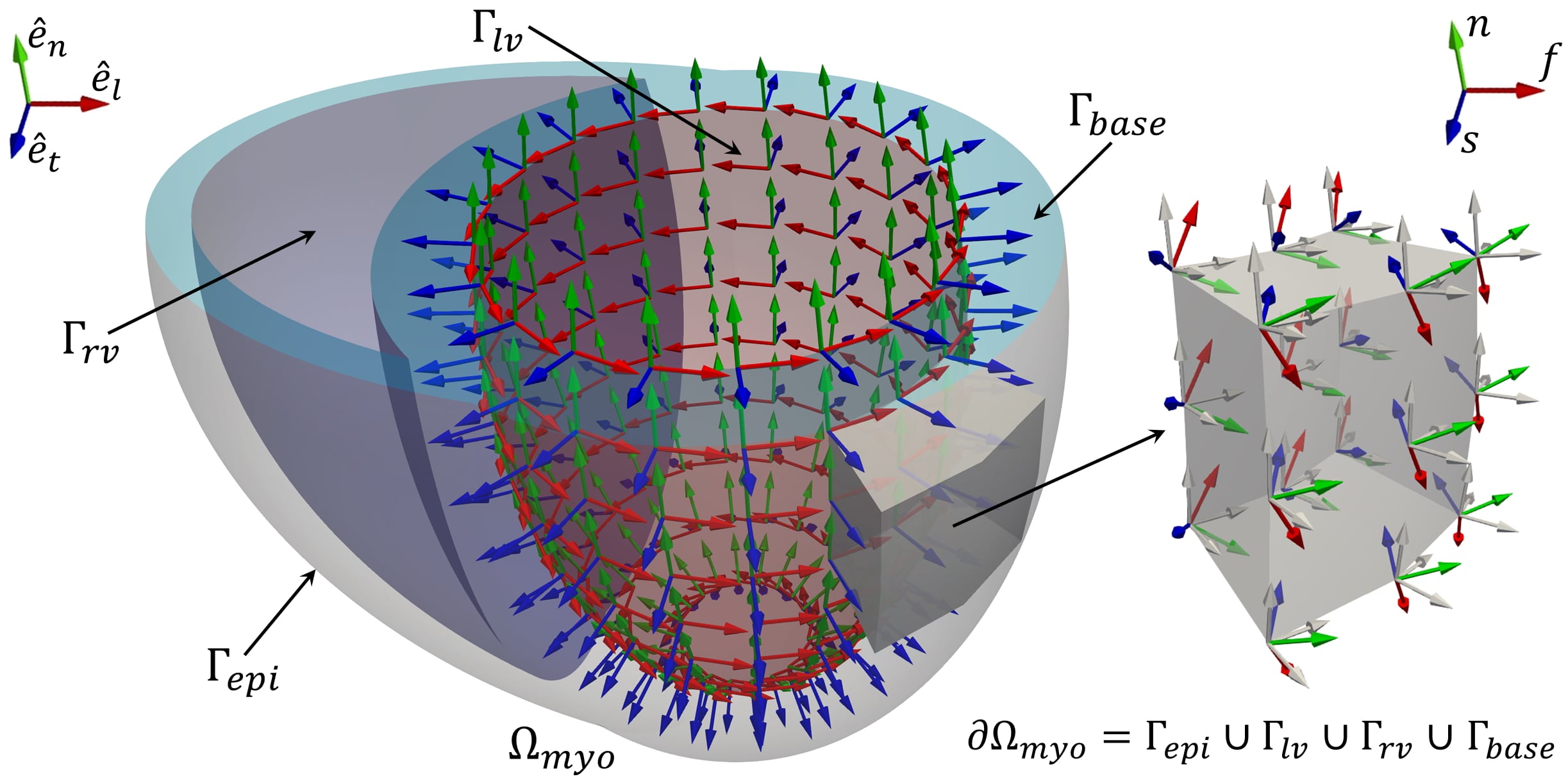}
	\caption{Left: Representation of the three directions employed by a LDRBM
		for an idealized biventricular domain $\Omega_{myo}$ whose border $\partial \Omega_{myo}$ is partitioned in the epicardium $\Gamma_{epi}$, left $\Gamma_{lv}$ and right $\Gamma_{rv}$ endocardium and the \textcolor{black}{ventricular} base $\Gamma_{base}$, $\partial \Omega_{myo}=\Gamma_{epi} \cup \Gamma_{lv} \cup \Gamma_{rv} \cup \Gamma_{base}$. For visualization purpose, only directions on the left endocardium $\Gamma_{lv}$ are represented.
		In blue: unit transmural direction, $\widehat{\boldsymbol{e}}_t$; In green: unit normal direction, $\widehat{\boldsymbol{e}}_n$; In red: unit longitudinal direction, $\widehat{\boldsymbol{e}}_l$. 
		Right: zoom on a slab of the left \textcolor{black}{ventricular} myocardium $\Omega_{myo}$ showing the three final myofibers orientations $\boldsymbol f,\,\boldsymbol s,\,\boldsymbol n$.}
	\label{fig:directions}
\end{figure}

We identify the following shared steps of the three ventricular LDRBMs which are hereby reported:
\begin{description}
	\item[1. Labelled mesh:] Provide a labelled mesh of the \textcolor{black}{ventricular} domain $\Omega_{myo}$ to define specific partitions of the boundary $\partial \Omega_{myo}$, see Figure \ref{fig:directions}; 
	\item[2. Transmural distance:] A {\sl transmural distance} is defined to compute the distance of the epicardium from endocardial surfaces;
	\item[3. Transmural direction:] The transmural distance gradient is used to build the unit {\sl transmural direction}~$\widehat{\boldsymbol{e}}_t$
	of the ventricles, see Figure \ref{fig:directions};
	\item[4. Normal direction:] An apico-basal direction (directed from the apex towards the \textcolor{black}{ventricular} base) is introduced and it is used to build the unit {\sl normal direction} $\widehat{\boldsymbol{e}}_n$, orthogonal to the transmural one, see Figure \ref{fig:directions};
	\item[5. Local coordinate system:] Build for each point of the \textcolor{black}{ventricular} domain an orthonormal local coordinate axial system composed by
	$\widehat{\boldsymbol{e}}_t,\,\widehat{\boldsymbol{e}}_n$ and the unit {\sl longitudinal direction} $\widehat{\boldsymbol{e}}_l$ (orthogonal to the previous ones), see Figure \ref{fig:directions}; 
	\item[6. Rotate axis:] Finally, properly rotate the reference frame with the purpose of defining the myofiber orientations: $\boldsymbol f$ the {\sl fiber direction}, $\boldsymbol n$ the {\sl cross-fiber direction} and $\boldsymbol s$ the {\sl sheet direction}, see Figure~\ref{fig:directions}(right). Rotations are chosen in order to match histology and DTI observations.
\end{description}
The ventricular LDRBM by Rossi et al. \cite{rossi2014thermodynamically, quarteroni2017integrated} (in what follows referred to as R-RBM) is a modified version of the algorithm studied in \cite{wong2014generating}  for generating fibers in left ventricular geometries \cite{rossi2014thermodynamically}, then extended to the biventricular case in \cite{quarteroni2017integrated}. R-RBM is based on the definition of the transmural direction. Bayer et al.~\cite{bayer2012novel} developed another LDRBM (B-RBM) for assigning myocardial fiber orientation introducing two major contributions. The first improvement is the definition of the apico-basal direction \cite{bayer2005laplace}. The second one consists of using the bi-direction spherical interpolation (\texttt{bislerp}) \cite{shoemake1985animating, kuipers1999quaternions} to manage the fiber orientations in order to guarantee a smooth and continuous change in the fiber field, particularly in the proximity of the septum and around the \textcolor{black}{ventricular} junctions \cite{bayer2012novel}. Both R-RBM and B-RBM introduce an artificial basal plane, located well below the cardiac valves, delimiting the \textcolor{black}{ventricular} regions. To overcome this restriction, Doste et al. \cite{doste2019rule} recently proposed a new LDRBM (D-RBM) which is able to build the fiber field in a full biventricular geometry without the need to cut it with a basal plane. D-RBM fiber orientations are generated to take into account specific ventricular regions, such as the inter-ventricular septum and outflow tracts (OT), following observations from histological studies \cite{doste2019rule}.

To characterize the three LDRBMs under review, it is useful to consider the following Laplace-Dirichlet problem
\begin{equation}
\label{laplace}
\begin{cases}
-\Delta \chi=0 &\qquad{\text{in }}\Omega_{myo},
\\
\chi = \chi_{a} &\qquad{\text{on }}\Gamma_{a},
\\
\chi = \chi_{b} &\qquad{\text{on }}\Gamma_{b}, 
\\
\nabla \chi \cdot \mathbf{n}=0 &\qquad{\text{on }}\Gamma_{n},
\end{cases}
\end{equation}  
for a generic unknown $\chi$ and suitable boundary data $\chi_{a},\,\chi_{b} \in \mathbb{R}$ set on generic partitions of the \textcolor{black}{ventricular} boundary $\Gamma_{a},\,\Gamma_{b},\,\Gamma_{n}$, with $\Gamma_{a}\cup\Gamma_{b}\cup\Gamma_{n}=\partial\Omega_{myo}$. The variable $\chi$ will assume different meanings depending on the step and LDRBM considered. Moreover, the values $\chi_{a},\,\chi_{b}$ are fixed in order to evaluate specific ventricular distances between boundary partitions $\Gamma_{a},\,\Gamma_{b}$.

\begin{figure}[t]
	\centering
	\includegraphics[width=1\textwidth]{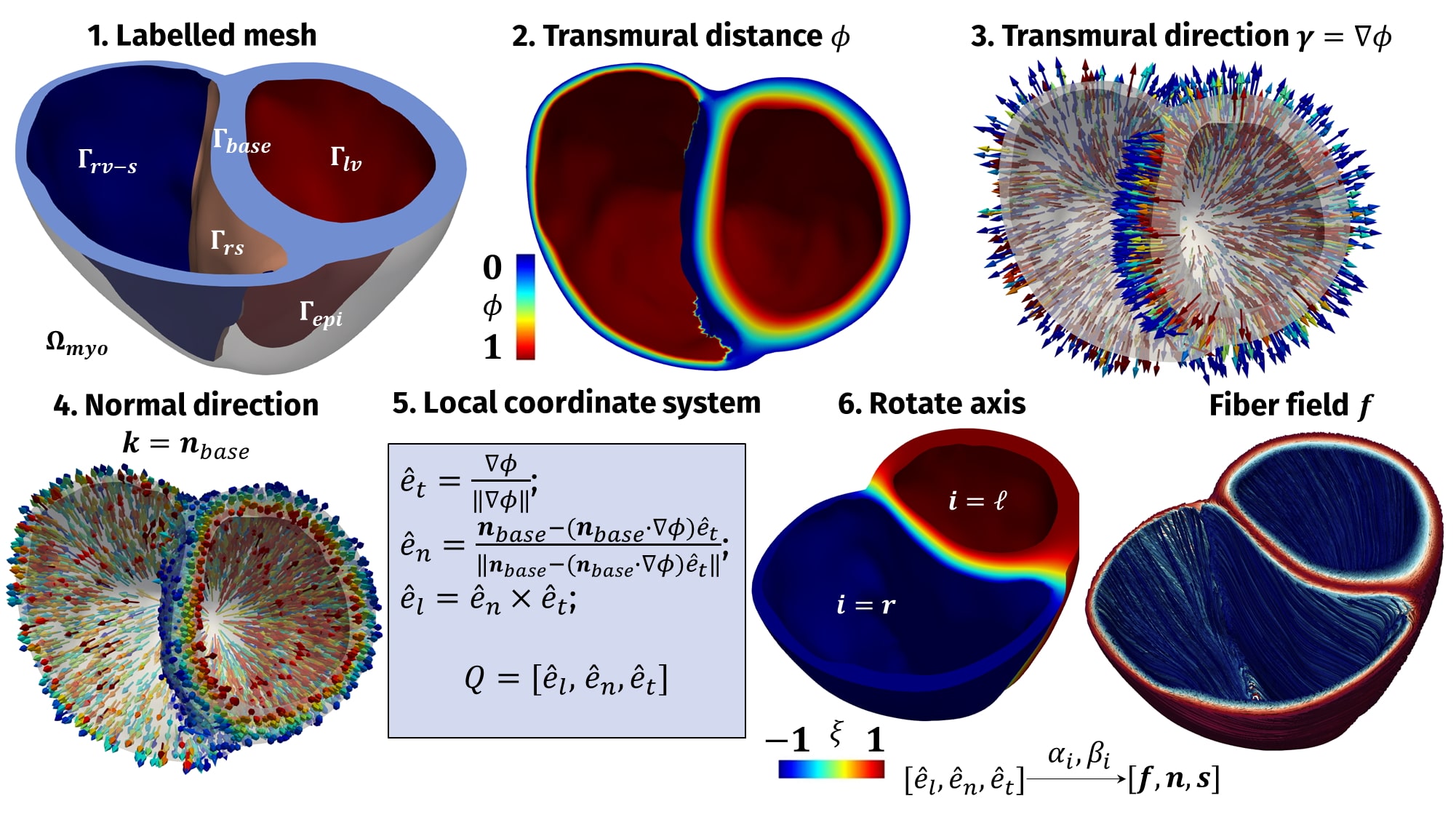}
	\caption{Schematic pipeline of R-RBM for a biventricular geometry with an artificial basal plane.}
	\label{fig:Rossi-LDRBM}
\end{figure}

\begin{figure}[t]
	\centering
	\includegraphics[width=1\textwidth]{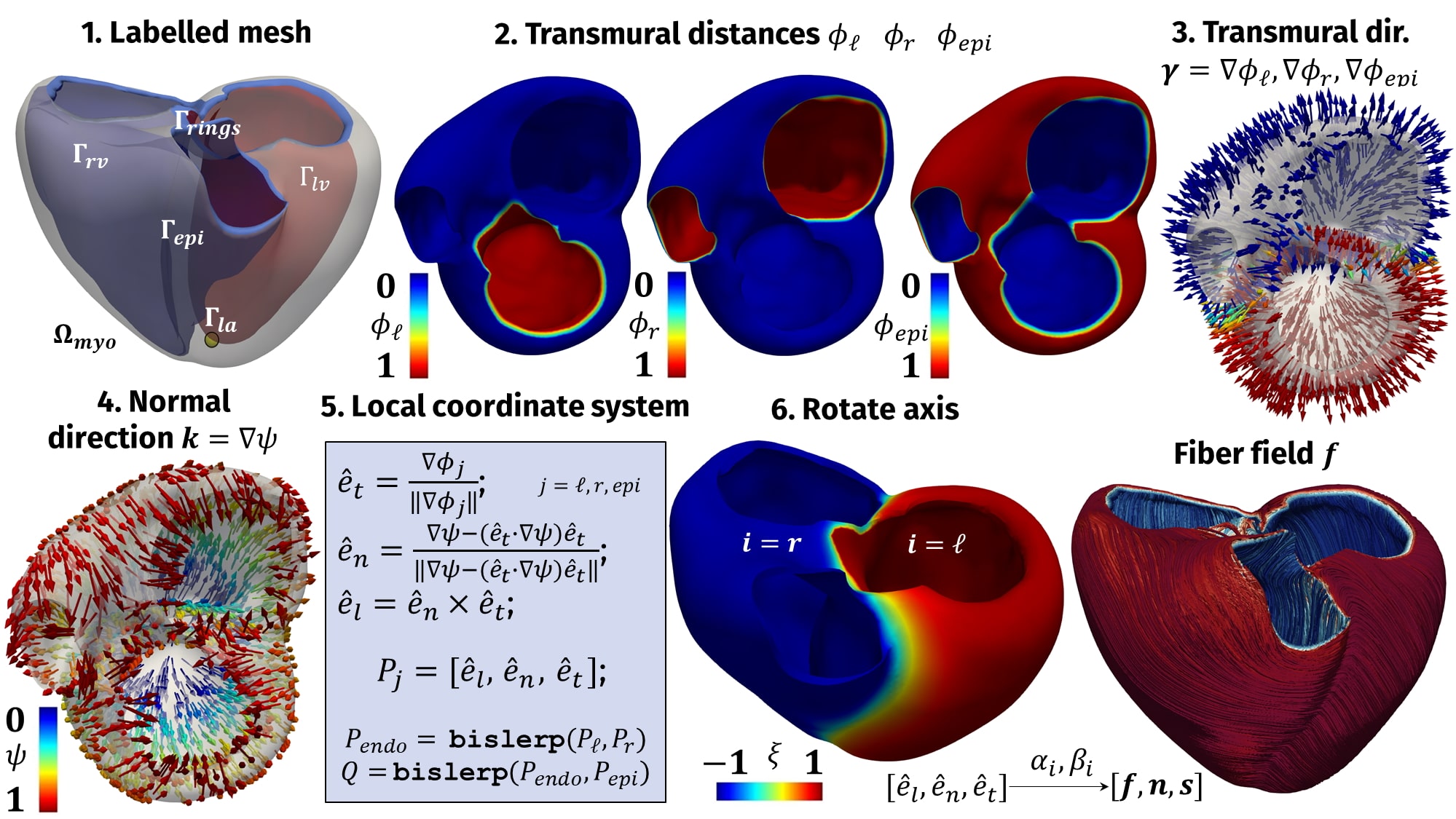}
	\caption{Schematic pipeline of B-RBM for a full biventricular geometry.}
	\label{fig:Bayer-LDRBM}
\end{figure}

\begin{figure}[t]
	\centering
	\includegraphics[width=1\textwidth]{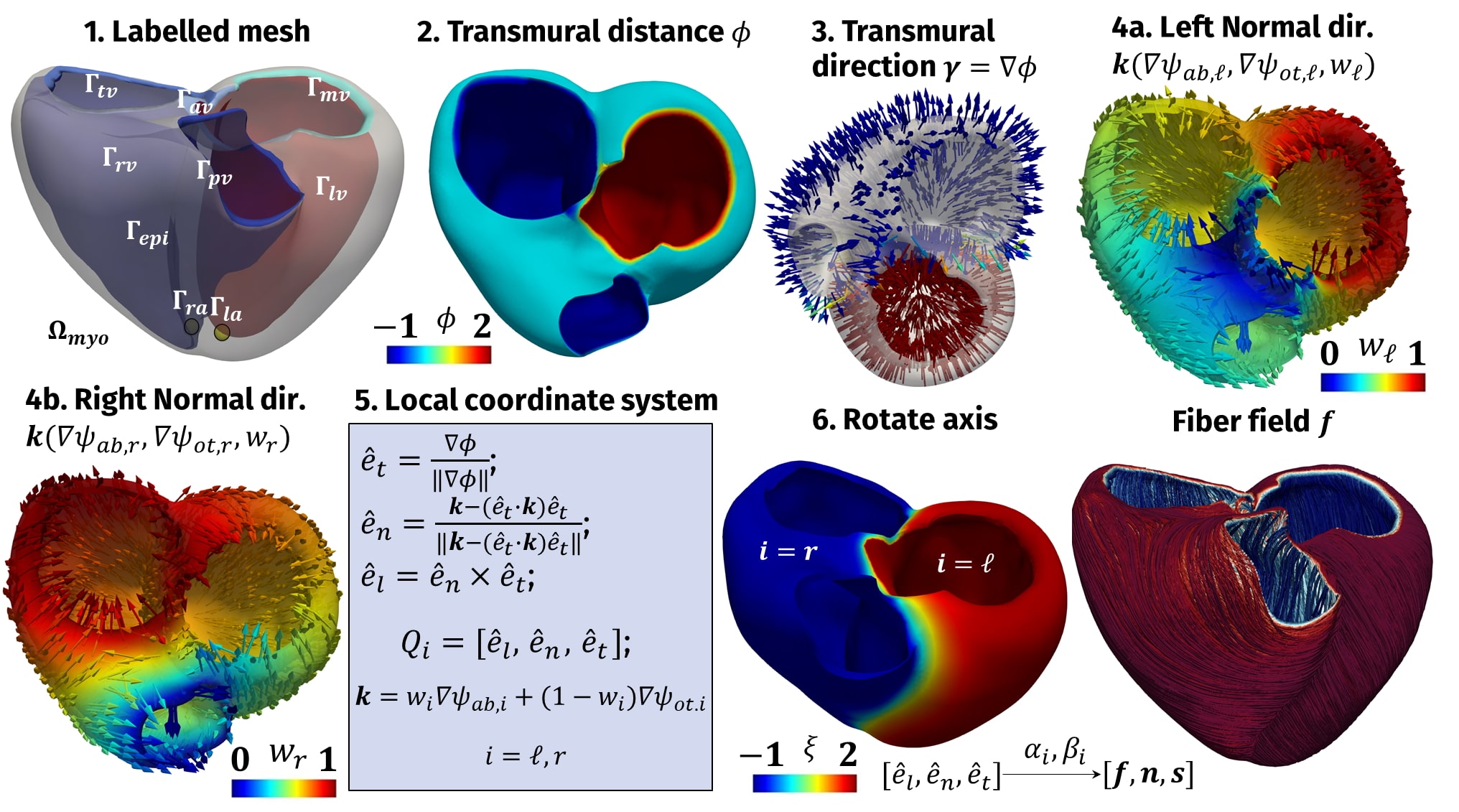}
	\caption{Schematic pipeline of D-RBM for a full biventricular geometry.}
	\label{fig:Doste-LDRBM}
\end{figure}

We detail in what follows the six points aforementioned. We refer to Figures \ref{fig:Rossi-LDRBM}, \ref{fig:Bayer-LDRBM} and \ref{fig:Doste-LDRBM}, showing a schematic representations of R-RBM, B-RBM, and D-RBM, respectively, for a biventricular domain $\Omega_{myo}$.

\paragraph{1. Labelled mesh}
Prescribing the \textcolor{black}{ventricular} boundary $\partial\Omega_{myo}$. All the three LDRBMs define the following boundaries:
\begin{equation*}
\begin{alignedat}{5}
\Gamma_{epi} &: {\text{ the epicardium}} ,
& \qquad
\Gamma_{base} &: {\text{ the ventricular base}},
\\
\Gamma_{lv} &: {\text{ the left endocardium}} ,
& \qquad
\Gamma_{rv} &: {\text{ the right endocardium}}.
\end{alignedat}
\end{equation*} 

Moreover, R-RBM subdivides the right \textcolor{black}{ventricular} endocardium $\Gamma_{rv}$ into the right septum $\Gamma_{rs}$ and the remaining part $\Gamma_{rv-s}$ such that $\Gamma_{rv}=\Gamma_{rs} \cup \Gamma_{rv-s}$, see step 1 in Figure \ref{fig:Rossi-LDRBM}. This subdivision is usually performed manually by the user, thus introducing some arbitrariness during the septum selection. For B-RBM and D-RBM the left \textcolor{black}{ventricular} apex $\Gamma_{la}$ is also introduced, whereas the right \textcolor{black}{ventricular} apex $\Gamma_{ra}$ for D-RBM solely (see step 1 in Figures \ref{fig:Bayer-LDRBM} and \ref{fig:Doste-LDRBM}). Furthermore, D-RBM requires boundary tags for the four valve rings: $\Gamma_{mv}$ (mitral valve), $\Gamma_{av}$ (aortic valve), $\Gamma_{tv}$ (tricuspid valve) and $\Gamma_{pv}$ (pulmonary valve), see step 1 in Figure \ref{fig:Doste-LDRBM}. It is also useful to define $\Gamma_{rings}=\Gamma_{lring} \cup \Gamma_{rring}$, with $\Gamma_{lring}=\Gamma_{mv} \cup \Gamma_{av}$ and $\Gamma_{rring}=\Gamma_{tv} \cup \Gamma_{pv}$. Notice that in B-RBM we considered the union of the four valve rings as the \textcolor{black}{ventricular} base $\Gamma_{base}=\Gamma_{rings}$. This allows the use of B-RBM also in the case of a full biventricular geometry, see step 1 in Figure \ref{fig:Bayer-LDRBM}. In summary, the three methods define the boundary $\partial\Omega_{myo}$ as follows (see step 1 in Figures~\ref{fig:Rossi-LDRBM},~\ref{fig:Bayer-LDRBM} and \ref{fig:Doste-LDRBM}; \textcolor{black}{see also the Appendix for further details about the tagging procedure}):    
\begin{equation*}
\begin{alignedat}{3}
& \text{R-RBM}: \partial\Omega_{myo} = \Gamma_{epi} \cup \Gamma_{lv} \cup \Gamma_{rs} \cup \Gamma_{rv-s} \cup \Gamma_{base},
\\
& \text{B-RBM}: \partial\Omega_{myo} = \Gamma_{epi} \cup \Gamma_{lv} \cup \Gamma_{rv} \cup \Gamma_{rings} \cup \Gamma_{la},
\\
& \text{D-RBM}: \partial\Omega_{myo} = \Gamma_{epi} \cup \Gamma_{lv} \cup \Gamma_{rv} \cup \Gamma_{mv} \cup \Gamma_{av} \cup \Gamma_{tv} \cup \Gamma_{pv} \cup \Gamma_{la} \cup \Gamma_{ra}.  
\end{alignedat}
\end{equation*}  

\begin{table}[t!]
	\label{table_phi}
	\begin{center}
		\begin{tabular}{ |c|c|c|c|c|c|c| } 
			\hline
			LDRBM type & Transmural distances & $\chi_{a}$ & $\Gamma_{a}$ & $\chi_{b}$ & $\Gamma_{b}$ & $\Gamma_{n}$  \\
			\hline
			\hline
			R-RBM & $\phi$ & 1 & $\Gamma_{lv} \cup \Gamma_{rv-s}$ & 0 & $\Gamma_{epi} \cup \Gamma_{rs}$ & $\Gamma_{base}$  \\
			\hline
			\multirow{3}{*}{B-RBM} & $\phi_{\ell}$  & 1 & $\Gamma_{lv}$ & 0 & $\Gamma_{epi} \cup \Gamma_{rv}$ & $\Gamma_{rings}$ \\ 
			& $\phi_{r}$  & 1 & $\Gamma_{rv}$ & 0 & $\Gamma_{epi} \cup \Gamma_{lv}$ & $\Gamma_{rings}$ \\ 
			& $\phi_{epi}$ & 1 & $\Gamma_{epi}$ & 0 & $\Gamma_{lv} \cup \Gamma_{rv}$ & $\Gamma_{rings}$ \\
			\hline 
			\multirow{2}{*}{D-RBM} & \multirow{2}{*}{$\phi$} & 2 & $\Gamma_{lv}$ & \multirow{2}{*}{0} & \multirow{2}{*}{$\Gamma_{epi}$} & \multirow{2}{*}{$\Gamma_{rings}$} \\
			& & -1 & $\Gamma_{rv}$ & & & \\
			\hline
		\end{tabular}
	\end{center}
	\caption{Transmural distance boundary conditions for R-RBM, B-RBM and D-RBM used in step 2.}
\end{table}
\paragraph{2. Transmural distance}
Definition of transmural distances (generally indicated with the letter $\phi$) obtained by solving Laplace-Dirichlet problems of the form \eqref{laplace}. In particular, for R-RBM, the transmural distance $\phi$ is found by solving \eqref{laplace} with $\chi_{a}=1$ on $\Gamma_{lv} \cup \Gamma_{rv-s}$, $\chi_{b}=0$ on 
$\Gamma_{epi} \cup \Gamma_{rs}$,  and $\Gamma_{n}=\Gamma_{base}$. For D-RBM,  $\phi$ is found by solving \eqref{laplace} with $\chi_{a}=2$ on $\Gamma_{lv}$, $\chi_{a}=-1$ on $\Gamma_{rv}$, $\chi_{b}=0$ on $\Gamma_{epi}$, and $\Gamma_{n}=\Gamma_{rings}$. 
B-RBM requires to solve three Laplace problems \eqref{laplace} in order to compute three different transmural distances $\phi_{\ell}$, $\phi_{r}$ and $\phi_{epi}$. We refer the reader to Table 1 for the specific choices in problem \eqref{laplace} made by the three methods. Notice that in D-RBM the boundary conditions $\chi_{a}$ are assigned in order to identify the two ventricles (positive and negative values for left and right ventricle, respectively) and to associate roughly two-thirds of the septum to the left ventricle and one-third to the right one \cite{doste2019rule} (see step 2 in Figures \ref{fig:Rossi-LDRBM}, \ref{fig:Bayer-LDRBM} and \ref{fig:Doste-LDRBM}).

\paragraph{3. Transmural direction}
After solving the Laplace problems for finding the transmural distances $\phi,\,\phi_l,\,\phi_r,\,\phi_{epi}$, their gradients define the transmural directions $\xx\gamma$ (see step 3 in Figures \ref{fig:Rossi-LDRBM}, \ref{fig:Bayer-LDRBM} and \ref{fig:Doste-LDRBM}). In particular, we have:
\begin{equation*}
\begin{alignedat}{3}
& \text{R-RBM}: \boldsymbol{\gamma} = \nabla\phi,
\\
& \text{B-RBM}: \boldsymbol{\gamma} = \nabla \phi_i,\quad i=\ell,r,epi,
\\
& \text{D-RBM}: \boldsymbol{\gamma} = \nabla\phi.  
\end{alignedat}
\end{equation*}

\paragraph{4. Normal direction}
Definition of the normal direction $\boldsymbol{k}$. In R-RBM, this is done using the vector $\mathbf{n}_{base}$, i.e. the outward normal to the ventricular base, that is $\boldsymbol{k}=\mathbf{n}_{base}$ (see step 4 in Figure \ref{fig:Rossi-LDRBM}). For the other two LDRBMs, further variables (identified by the letter $\psi$) are introduced: they are found by solving the Laplace problem in the form \eqref{laplace} with suitable boundary conditions. Specifically, for B-RBM the vector $\boldsymbol{k}$ is the gradient of the solution $\psi$ ($\boldsymbol{k}=\nabla \psi$) obtained with $\chi_{a}=1$ on $\Gamma_{rings}$, $\chi_{b}=0$ on $\Gamma_{la}$, and $\Gamma_{n}=\Gamma_{epi} \cup \Gamma_{lv} \cup \Gamma_{rv}$, see step 4 in Figure \ref{fig:Bayer-LDRBM}. In D-RBM, instead, two normal directions are introduced, one for each ventricle:
\begin{equation}
\label{doste_normal}
\boldsymbol{k} = w_i\nabla \psi_{ab,i} + (1-w_i)\nabla \psi_{ot,i}, \qquad i=\ell,r,
\end{equation}
where $i=\ell,r$ refer to the left and right ventricle, respectively, so that the normal direction is a weighted sum of apico-basal ($\nabla \psi_{ab,i}$) and apico-outflow-tract ($\nabla \psi_{ot,i}$) directions, obtained using an inter-ventricular interpolation function $w_i$  \cite{doste2019rule}; the latter are given again by solutions of problems like \eqref{laplace} (see step 4a and 4b in Figure \ref{fig:Doste-LDRBM}). In particular, $\psi_{ab,\ell}$, $\psi_{ab,r}$, $\psi_{ot,\ell}$, $\psi_{ot,r}$, $w_{\ell}$, $w_r$ are found by solving \eqref{laplace} with $\chi_{a}=1$ on $\Gamma_{a}$ and $\chi_{b}=0$ on $\Gamma_{b}$, where $\Gamma_{a}$ and $\Gamma_{b}$ are boundary subsets listed in Table 2.
Summing up, the different methods compute the normal direction $\boldsymbol{k}$ as follows (see step 4 in Figures \ref{fig:Rossi-LDRBM}, \ref{fig:Bayer-LDRBM} and \ref{fig:Doste-LDRBM}):
\begin{equation*}
\begin{alignedat}{3}
& \text{R-RBM}: \boldsymbol{k} = \mathbf{n}_{base},
\\
& \text{B-RBM}: \boldsymbol{k} = \nabla \psi,
\\
& \text{D-RBM}: \boldsymbol{k} = w_i\nabla \psi_{ab,i} + (1-w_i)\nabla \psi_{ot,i}, \qquad i=\ell,r.  
\end{alignedat}
\end{equation*}
\begin{table}[t]
	\label{table_Doste}
	\begin{center}
		\begin{tabular}{ |c|c|c|c| } 
			\hline
			Scalar potentials & $\Gamma_{a}$ & $\Gamma_{b}$ & $\Gamma_{n}$  \\
			\hline
			\hline
			$\psi_{ab,\ell}$ & $\Gamma_{mv}$ & $\Gamma_{la}$ & $\Gamma_{epi} \cup \Gamma_{lv} \cup \Gamma_{rv} \cup \Gamma_{av} \cup \Gamma_{rring}$ \\
			\hline
			$\psi_{ab,r}$ & $\Gamma_{tv}$ & $\Gamma_{ra}$ & $\Gamma_{epi} \cup \Gamma_{lv} \cup \Gamma_{rv} \cup \Gamma_{pv} \cup \Gamma_{lring}$ \\
			\hline
			$\psi_{ot,\ell}$ & $\Gamma_{av}$ & $\Gamma_{la}$ & $\Gamma_{epi} \cup \Gamma_{lv} \cup \Gamma_{rv} \cup \Gamma_{mv} \cup \Gamma_{rring}$ \\
			\hline
			$\psi_{ot,r}$ & $\Gamma_{pv}$ & $\Gamma_{ra}$ & $\Gamma_{epi} \cup \Gamma_{lv} \cup \Gamma_{rv} \cup \Gamma_{tv} \cup \Gamma_{lring}$ \\
			\hline			
			$w_{\ell}$ & $\Gamma_{mv} \cup \Gamma_{la}$ & $\Gamma_{av}$ & $\Gamma_{epi} \cup \Gamma_{lv} \cup \Gamma_{rv} \cup \Gamma_{rring}$ \\
			\hline		
			$w_r$ & $\Gamma_{tv} \cup \Gamma_{ra}$ & $\Gamma_{pv}$ & $\Gamma_{epi} \cup \Gamma_{lv} \cup \Gamma_{rv} \cup \Gamma_{lring}$ \\
			\hline					
		\end{tabular}
	\end{center}
	\caption{Scalar potentials used in D-RBM to build the normal direction.}
\end{table}
\paragraph{5. Local coordinate system}
Building an orthonormal local coordinate system (defined by letter $Q$) at each point of the domain $\Omega_{myo}$. All the three methods make use of the following function \texttt{axis}:
\begin{equation}
\label{axis}
P=[\widehat{\boldsymbol{e}}_l, \widehat{\boldsymbol{e}}_n, \widehat{\boldsymbol{e}}_t]=\texttt{axis}(\boldsymbol{k}, \xx\gamma)=
\begin{cases}
\widehat{\boldsymbol{e}}_t &= \frac{\xx\gamma}{\left\lVert \xx\gamma \right\rVert}, 
\\
\widehat{\boldsymbol{e}}_n &= \frac{\boldsymbol{k} - (\boldsymbol{k} \cdot \widehat{\boldsymbol{e}}_t )\widehat{\boldsymbol{e}}_t}{\left\lVert \boldsymbol{k} - (\boldsymbol{k} \cdot \widehat{\boldsymbol{e}}_t )\widehat{\boldsymbol{e}}_t \right\rVert}, 
\\
\widehat{\boldsymbol{e}}_l &= \widehat{\boldsymbol{e}}_n \times \widehat{\boldsymbol{e}}_t,  
\end{cases}
\end{equation}
which takes as input a normal direction $\boldsymbol{k}$ and a transmural direction $\xx\gamma$ and returns 
the orthonormal system $P$ whose columns are the three orthonormal directions $\widehat{\boldsymbol{e}}_l$, $\widehat{\boldsymbol{e}}_n$, $\widehat{\boldsymbol{e}}_t$ which represent the longitudinal, the normal and the transmural unit directions, respectively. 
For R-RBM we have $Q=\texttt{axis}(\boldsymbol{k}, \nabla\phi)$. 
For B-RBM three orthonormal coordinate systems are introduced, that is  $P_{\ell}=\texttt{axis}(\boldsymbol{k}, \nabla \phi_{\ell})$, $P_r=\texttt{axis}(\boldsymbol{k}, \nabla \phi_r)$ and $P_{epi}=\texttt{axis}(\boldsymbol{k}, \nabla \phi_{epi})$, which are then interpolated through the function \texttt{bislerp}
to obtain a continuous orthonormal coordinate system within the whole myocardium (see \cite{bayer2012novel} for further details). Hence, B-RBM performs the following steps 
to obtain the final orthonormal coordinate system $Q$ (see step 5 in Figure \ref{fig:Bayer-LDRBM}):
\begin{equation*}
\begin{alignedat}{2}
P_{endo}=\texttt{bislerp}(P_{\ell},P_r),
\\
Q=\texttt{bislerp}(P_{endo},P_{epi}).
\\ 
\end{alignedat}
\end{equation*}
D-RBM, instead, defines two different coordinate systems for left and right ventricles as a consequence of the normal directions defined in \eqref{doste_normal} (see step 5 in Figure \ref{fig:Doste-LDRBM}):
\begin{equation*}
Q_{i}=\texttt{axis}(w_i\nabla \psi_{ab,i} + (1-w_i)\nabla \psi_{ot,i}, \nabla \phi), \quad i=\ell,r.
\end{equation*}

\paragraph{6. Rotate axis}
The orthonormal coordinate system, defined for each point of the myocardium at the previous step, should
be aligned in order to match histological knowledge about fiber and sheet orientations. 
To this aim, the three LDRBMs introduce a rotation of $\ww{\xx e}_l,\,\ww{\xx e}_n,\,\ww{\xx e}_t$ by means of suitable angles: the longitudinal direction $\widehat{\boldsymbol{e}}_l$ rotates counter-clockwise around $\widehat{\boldsymbol{e}}_t$ by an angle $\alpha_i$, whereas the transmural direction $\widehat{\boldsymbol{e}}_t$ is rotated counter-clockwise around $\widehat{\boldsymbol{e}}_l$ by an angle $\beta_i$,
where $i=\ell,r$ depend on the left or right ventricle the point belongs to. 
Indeed, it is known that in the left and right ventricles the fiber orientations feature a change in direction at the inter-ventricular septum \cite{kocica2006helical}.
In order to obtain a
local orthonormal coordinate system, direction $\widehat{\boldsymbol{e}}_n$ is rotated accordingly.

These rotations produce a map from the original coordinate system to a new coordinate system $[\boldsymbol f, \boldsymbol n, \boldsymbol s]$:
\begin{equation*}
[\widehat{\boldsymbol{e}}_l, \widehat{\boldsymbol{e}}_n, \widehat{\boldsymbol{e}}_t] \xrightarrow{\alpha_{i}, \beta_{i}} [\boldsymbol f, \boldsymbol n, \boldsymbol s],\quad i=\ell,\,r,
\end{equation*}
where $\boldsymbol f$ is the fiber direction, $\boldsymbol n$ is the cross-fiber direction and $\boldsymbol s$ is the sheet direction. 

For all the three methods the rotation angles $\alpha_i=\alpha_i(d_i)$ and $\beta_i=\beta_i(d_i)$ 
are functions of the position within the myocardium, in particular of the {\sl
	transmural normalized distance} $d_i\in [0,1],\,i=\ell,r,$ defined as: 
\begin{equation*}
\begin{alignedat}{3}
& \text{R-RBM}: d_{\ell} = d_{r} = \phi ,
\\
& \text{B-RBM}: d_{\ell} = \phi_{\ell} \quad d_{r} = \phi_r,
\\
& \text{D-RBM}: d_{\ell} = \phi/2 \quad d_{r} = |\phi|.  
\end{alignedat}
\end{equation*} 
Accordingly, the rotation angles are written by means of the
following linear relationships:
\begin{equation}
\label{angle}
\alpha_i(d_i) = \alpha_{epi,i}(1-d_i)+\alpha_{endo,i}d_i,  \qquad  
\beta_i(d_i) = \beta_{epi,i}(1-d_i)+\beta_{endo,i}d_i,
\qquad i=\ell,r,           
\end{equation}
where $\alpha_{endo,\ell}$, $\alpha_{epi,\ell}$, $\alpha_{endo,r}$, $\alpha_{epi,r}$, $\beta_{endo,\ell}$, $\beta_{epi,\ell}$, $\beta_{endo,r}$ and $\beta_{epi,r}$ are suitable rotation angles on 
the epicardium and endocardium chosen in order to match histological observations. For example, classical values found in the literature are $\alpha_{epi,\ell}=+60^o,\,
\alpha_{endo,\ell}=-60^o,\,\alpha_{epi,r}=-25^o,\,\alpha_{endo,r}=+90^o$ \cite{lombaert2012human,greenbaum1981left,anderson2009three,ho2006anatomy,sanchez2015anatomical}.

In order to differentiate between the left and right ventricles and to apply the correct angles, we propose here to use the solution of an additional Laplace 
problem \eqref{laplace} in the unknown $\chi=\xi$ with $\chi_a=1$ on $\Gamma_{lv}$, $\chi_b=-1$ on $\Gamma_{rv}$, and $\Gamma_{n}=\Gamma_{base} \cup \Gamma_{epi}$\footnote{Let us observe that, for B-RBM $\Gamma_{base}=\Gamma_{rings}$ in the case of a full biventricular geometry. Moreover, for D-RBM solely $\chi_a=2$ in order to be compliant with the transmural distance.}. In particular, positive values of $\xi$ identify the left ventricle, whereas negative values the right one \cite{bayer2018universal}. This new feature enables to perform different rotations for left and right ventricles (see steps 6 in Figures \ref{fig:Rossi-LDRBM}, \ref{fig:Bayer-LDRBM} and \ref{fig:Doste-LDRBM}) that is crucial in order to generate realistic fiber fields. An alternative method has been proposed in \cite{doste2019rule} but only for D-RBM. It is worth mentioning that the original R-RBM \cite{rossi2014thermodynamically,quarteroni2017integrated} introduces a rotation to obtain the fiber field $\xx f$ only. Here we propose an extension in order to define also $\boldsymbol n$ and $\boldsymbol s$.

Further, B-RBM exploits two other functions representing the rotation angles within the septum:
\begin{equation*}
\alpha_{s}(d_i) = \alpha_{endo,i}(1-2d_i),
\qquad
\beta_{s}(d_i) = \beta_{endo,i}(1-2d_i),
\qquad
i=\ell,r,
\end{equation*} 
whereas with similar expressions, D-RBM introduces also the possibility to set specific fiber and sheet angles rotation in the OT regions (see \cite{bayer2012novel} and \cite{doste2019rule} for further details).

\medskip

\noindent We conclude pointing out that B-RBM and D-RBM can be applied to the full biventricular geometry and to the based biventricular case (that is obtained with an artificial basal plane well below the cardiac valves). Indeed, in the based biventricular geometry the whole procedure for B-RBM and D-RBM remains the same as long as the ring tags are replaced by the base tag, $\Gamma_{rings}=\Gamma_{base}$. On the contrary, R-RBM is less suitable for a full biventricular case because it is not able to strictly identify the normal direction $\boldsymbol{k}$ as the outward normal to the ventricular rings. Besides, the right septum definition $\Gamma_{rs}$ can be arbitrary for a full biventricular geometry.


\section{A new rule-based method for atrial fibers generation}\label{sec:atrial}

Atrial fibers architecture is very different from that of the ventricles where myofibers are aligned in a regular pattern \cite{streeter1969fiber}. Indeed, fibers in the atria are organized in individual bundles running along different directions throughout the wall chambers. Preferred orientation of myofibers in the human atria is characterized by multiple overlapping structures, which promote the formation of separate attached bundles~\cite{dossel2012computational}. This architecture has a strong influence in the electric signal propagation inside the atrial muscle \cite{krueger2011modeling,ferrer2015detailed,boineau1988demonstration,betts2002three,de2002electroanatomic,kruger2013personalized,maesen2013rearrangement}. 

Over the past years many histo-anatomical studies investigated the fibers arrangement of the atria revealing a very complex texture musculature \cite{papez1920heart, thomas1959muscular, ho2002atrial, ho2009importance, ho2012left, aslanidi2012application, sanchez2013standardized, hansen2017fibrosis}. Nevertheless, there is a paucity of imaging data on atrial fibers orientation with respect to the ventricles, mainly due to imaging difficulties in capturing the thin atrial walls \cite{dossel2012computational}. Only recently, ex vivo atrial fibers have been analysed owing to submillimeter {\sl Diffusion Tensor MRI} imaging \cite{pashakhanloo2016myofiber,zhao2015integration,zhao2017three}.   

In computational models of cardiac EP, atrial fiber orientations have been generated in specific regions either manually \cite{harrild2000computer,vigmond2001reentry,jacquemet2003study,seemann2006heterogeneous} or using a variety of semi-automatic rule-based approaches \cite{krueger2011modeling,tobon2013three,fastl2018personalized,ferrer2015detailed,plank2008evaluating,krueger2010patient,labarthe2012semi,rocher2019highly}. Recently, atlas-based methods, in which fiber directions of a reference atrial geometry are warped on a target geometry, have been introduced \cite{hoermann2019automatic,roney2020constructing,satriano2013feature,mcdowell2015virtual,roney2019universal}. All the former procedures require manual intervention introducing, for example, various distinct landmarks, seed-points and a network of auxiliaries lines \cite{krueger2011modeling,fastl2018personalized}. Hence, a processing pipeline for generating atrial fibers field still remains a knotty procedure~\cite{fastl2018personalized,fastl2016personalized}. 

In this section we propose for the first time a LDRBM for the generation of atrial myofibers, which is able to qualitatively reproduce all the important features,  such as fiber bundles, needed to provide a realistic atrial musculature architecture. Our newly developed method is inspired by \cite{roney2019universal} where Laplace
problems are introduced to map variables between two geometries and by LDRBMs,  purposely built for the ventricles, presented in Section \ref{sec:ventr} \cite{bayer2012novel,rossi2014thermodynamically,doste2019rule}. 
The extension of the latter is not straightforward due to the presence
of bundles which run in different directions. For this reason, our atrial LDRBM combines the gradient of several harmonic functions to represent the fiber bundles. 

\begin{figure}[ht!]
	\centering
	\includegraphics[width=0.88\textwidth]{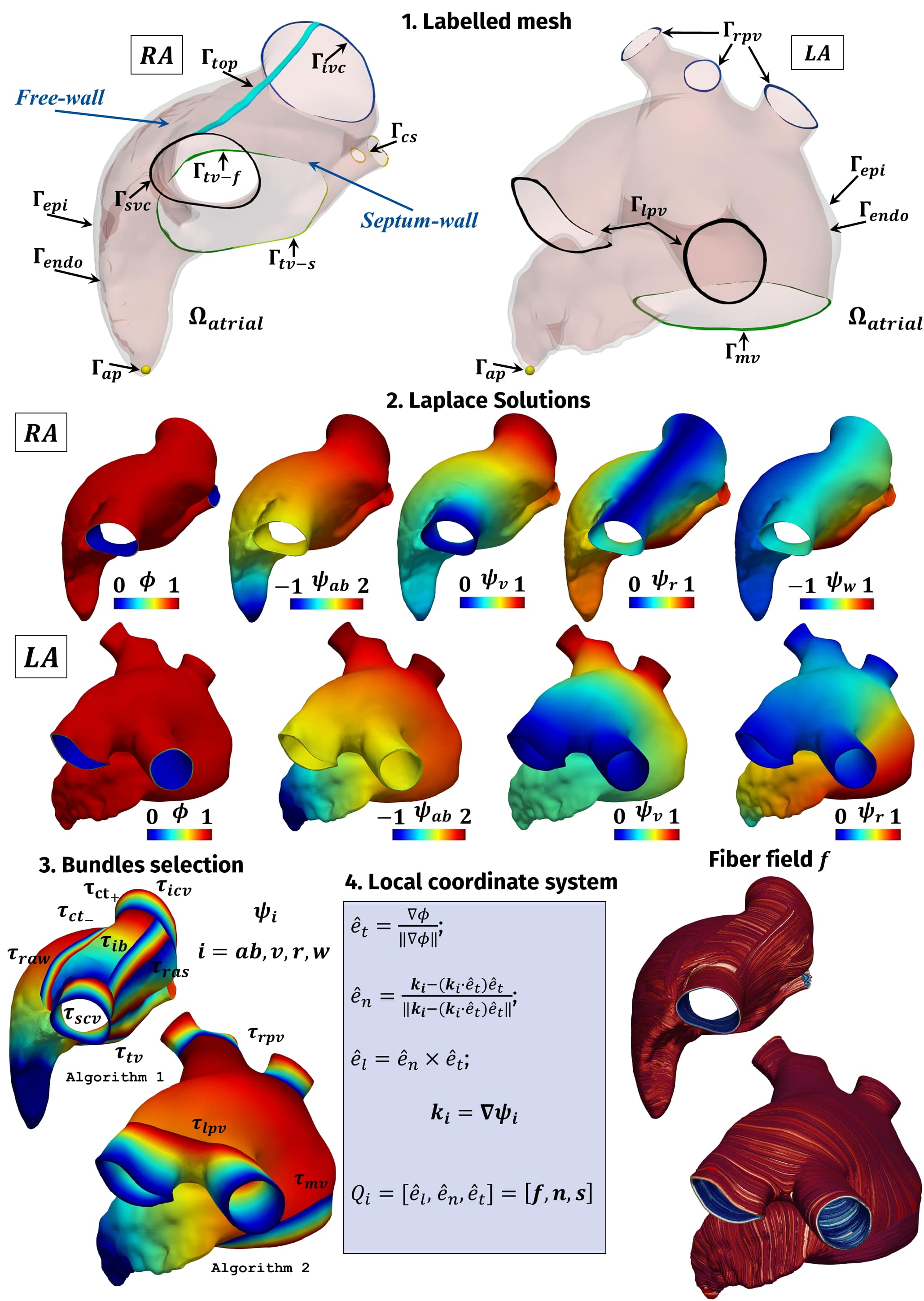}
	\caption{Schematic pipeline of the atrial LDRBM in a realistic right atrium (RA) and left atrium (LA) geometries.}
	\label{fig:LARA-LDRBM}
\end{figure}

In what follows we detail the four steps of the proposed atrial LDRBM.
We refer to Figure \ref{fig:LARA-LDRBM} for a schematic representation of the method in a realistic scenario.
\paragraph{1. Labelled mesh}
Label the mesh of the atrial computational domain $\Omega_{atrial}$ to define the boundary partitions of $\partial\Omega_{atrial}$. For both the left atrium (LA) and the right atrium (RA) the following boundaries are defined (see step 1 in Figure \ref{fig:LARA-LDRBM}):
\begin{equation*}
\begin{alignedat}{5}
\Gamma_{endo} &: {\text{ the endocardium}} ,
& \qquad
\Gamma_{epi} &: {\text{ the epicardium}} ,
& \qquad
\Gamma_{ap} &: {\text{ the apex appendage}}.
\end{alignedat}
\end{equation*}
\textcolor{black}{Moreover, for LA we introduce boundary tags of the mitral valve ring $\Gamma_{mv}$, the left and right pulmonary vein rings $\Gamma_{lpv}$, $\Gamma_{rpv}$. In the RA, the tags of the inferior and superior caval vein rings $\Gamma_{icv}$, $\Gamma_{scv}$, the coronary sinus ring $\Gamma_{cs}$ and the tricuspid valve ring $\Gamma_{tv}$ are introduced. In particular, $\Gamma_{tv}$ is equally subdivided in one part facing the atrial septum $\Gamma_{tv-s}$ and another one related to the free wall $\Gamma_{tv-f}$, such that $\Gamma_{tv}=\Gamma_{tv-s} \cup~\Gamma_{tv-f}$ (see step 1 in Figure \ref{fig:LARA-LDRBM}).} Furthermore, the RA encloses also the boundary tag \textcolor{black}{connecting the top upper region of the inferior and superior caval vein rings $\Gamma_{top}$} that is further divided in two sections: a part lying on the epicardium \textcolor{black}{$\Gamma_{top-epi}$} and another one on the endocardium \textcolor{black}{$\Gamma_{top-endo}$}, such that $\Gamma_{top}=\Gamma_{top-epi} \cup \Gamma_{top-endo}$ (see step 1 in Figure \ref{fig:LARA-LDRBM}; \textcolor{black}{see also the Appendix for further details about the tagging procedure}). In summary, the boundaries $\partial\Omega_{atrial}$ for the RA and LA are, respectively:  
\begin{equation*}
\begin{alignedat}{2}
& \text{RA}: \partial\Omega_{atrial} = \Gamma_{epi} \cup \Gamma_{endo} \cup \Gamma_{ap} \cup \Gamma_{icv} \cup \Gamma_{scv} \cup \Gamma_{cs} \cup \textcolor{black}{\Gamma_{tv-s} \cup \Gamma_{tv-f} \cup \Gamma_{top-epi} \cup \Gamma_{top-endo}};
\\
& \text{LA}: \partial\Omega_{atrial} = \Gamma_{epi} \cup \Gamma_{endo} \cup \Gamma_{ap} \cup \Gamma_{lpv} \cup \Gamma_{rpv} \cup \textcolor{black}{\Gamma_{mv}}.  
\end{alignedat}
\end{equation*}  
\begin{table}[t]
	\begin{center}
		\begin{tabular}{ |c|c|c|c|c|c|c| } 
			\hline
			Type & $\chi$ & $\chi_{a}$ & $\Gamma_{a}$ & $\chi_{b}$ & $\Gamma_{b}$ & $\Gamma_{n}$  \\
			\hline
			\hline
			\multirow{4}{*}{LA} & $\phi$ & 1 & $\Gamma_{epi}$ & 0 & $\Gamma_{endo}$ & \textcolor{black}{$\Gamma_{mv}$} $\cup \Gamma_{lpv} \cup \Gamma_{rpv}$ \\\cline{2-7} 
			& \multirow{2}{*}{$\psi_{ab}$} & 2 & $\Gamma_{rpv}$ & 1 & \textcolor{black}{$\Gamma_{mv}$} & \multirow{2}{*}{$\Gamma_{epi} \cup \Gamma_{endo}$} \\
			& & 0 & $\Gamma_{lpv}$ & -1 & $\Gamma_{ap}$ & \\\cline{2-7} 
			& $\psi_v$ & 1 & $\Gamma_{rpv}$ & 0 & $\Gamma_{lpv}$ & \textcolor{black}{$\Gamma_{mv}$} $\cup \Gamma_{epi} \cup \Gamma_{endo}$ \\\cline{2-7} 
			& $\psi_{r}$ & 1 & \textcolor{black}{$\Gamma_{mv}$} & 0 & $\Gamma_{lpv} \cup \Gamma_{rpv} \cup \Gamma_{ap}$ & $\Gamma_{epi} \cup \Gamma_{endo}$ \\                       
			\hline				
			\hline
			\multirow{5}{*}{RA} & $\phi$ & 1 & $\Gamma_{epi} \cup$ \textcolor{black}{$\Gamma_{top-epi}$} & 0 & $\Gamma_{endo} \cup$ \textcolor{black}{$\Gamma_{top-endo}$} & \textcolor{black}{$\Gamma_{tv}$} $\cup \Gamma_{icv} \cup \Gamma_{scv} \cup \Gamma_{cs}$ \\\cline{2-7} 
			& \multirow{2}{*}{$\psi_{ab}$} & 2 & $\Gamma_{icv}$ & 1 & \textcolor{black}{$\Gamma_{tv}$} & \multirow{2}{*}{$\Gamma_{epi} \cup \Gamma_{endo} \cup \textcolor{black}{\Gamma_{top}} \cup \Gamma_{cs}$} \\
			& & 0 & $\Gamma_{scv}$ & -1 & $\Gamma_{ap}$ & \\\cline{2-7} 
			& $\psi_v$ & 1 & $\Gamma_{icv}$ & 0 & $\Gamma_{scv} \cup \Gamma_{ap}$ & \textcolor{black}{$\Gamma_{mv}$} $\cup  \Gamma_{epi} \cup \Gamma_{endo} \cup \textcolor{black}{\Gamma_{top}} \cup \Gamma_{cs}$ \\\cline{2-7}
			& $\psi_r$ & 1 & \textcolor{black}{$\Gamma_{tv}$} & 0 & $\textcolor{black}{\Gamma_{top}}$ & $\Gamma_{epi} \cup \Gamma_{endo} \cup \Gamma_{icv} \cup \Gamma_{scv} \cup \Gamma_{cs}$ \\\cline{2-7}
			& $\psi_w$ & 1 & \textcolor{black}{$\Gamma_{tv-s}$} & -1 & \textcolor{black}{$\Gamma_{tv-f}$} & $\Gamma_{epi} \cup \Gamma_{endo} \cup \textcolor{black}{\Gamma_{top}} \cup \Gamma_{icv} \cup \Gamma_{scv} \cup \Gamma_{cs}$ \\                        
			\hline					
		\end{tabular}
	\end{center}
	\caption{Boundary data chosen in the Laplace problem \eqref{laplace} for the transmural distance $\phi$ and the \textcolor{black}{intra-atrial} distances $\psi_i$ ($i=ab,v,r,w$) in the left (LA) and right atrium (RA).}
	\label{table_atria}
\end{table}
\paragraph{2. Laplace solutions}
Definition of the transmural distance $\phi$ and several \textcolor{black}{intra-atrial} distances $\psi_i$, obtained by solving a Laplace-Dirichlet problem in the form of \eqref{laplace} with proper Dirichlet boundary conditions on the atrial boundaries, see step 2 in Figure \ref{fig:LARA-LDRBM}. Refer to Table \ref{table_atria} for the specific choices in problem \eqref{laplace} made for the RA and LA. In particular, $\psi_{ab}$ is the solution of a Laplace problem \eqref{laplace} with three different boundary data prescribed on the atrial appendage $\Gamma_{ap}$, the rings of the caval veins ($\Gamma_{scv}, \Gamma_{icv}$) \textcolor{black}{and the tricuspid-valve ring $\Gamma_{tv}$} for the RA, the pulmonary veins ($\Gamma_{lpv}, \Gamma_{rpv}$) and \textcolor{black}{mitral-valve ring $\Gamma_{mv}$} for LA; $\psi_v$ represents the distance between the caval veins for the RA and among the pulmonary veins for LA; $\psi_r$ stands for the distance between \textcolor{black}{the tricuspid-valve ring $\Gamma_{tv}$} and $\textcolor{black}{\Gamma_{top}}$ (RA) and between \textcolor{black}{the mitral-valve ring $\Gamma_{mv}$} and the union of the pulmonary veins rings $\Gamma_{lpv}\cup\Gamma_{rpv}$ (LA). Moreover, for the RA $\psi_w$ is the distance between the \textcolor{black}{tricuspid-valve ring of the free ($\Gamma_{tv-f}$) and the septum ($\Gamma_{tv-s}$) walls}. See step 2 in Figure \ref{fig:LARA-LDRBM}.    
\begin{figure}[t!]
	\centering
	\includegraphics[width=1\textwidth]{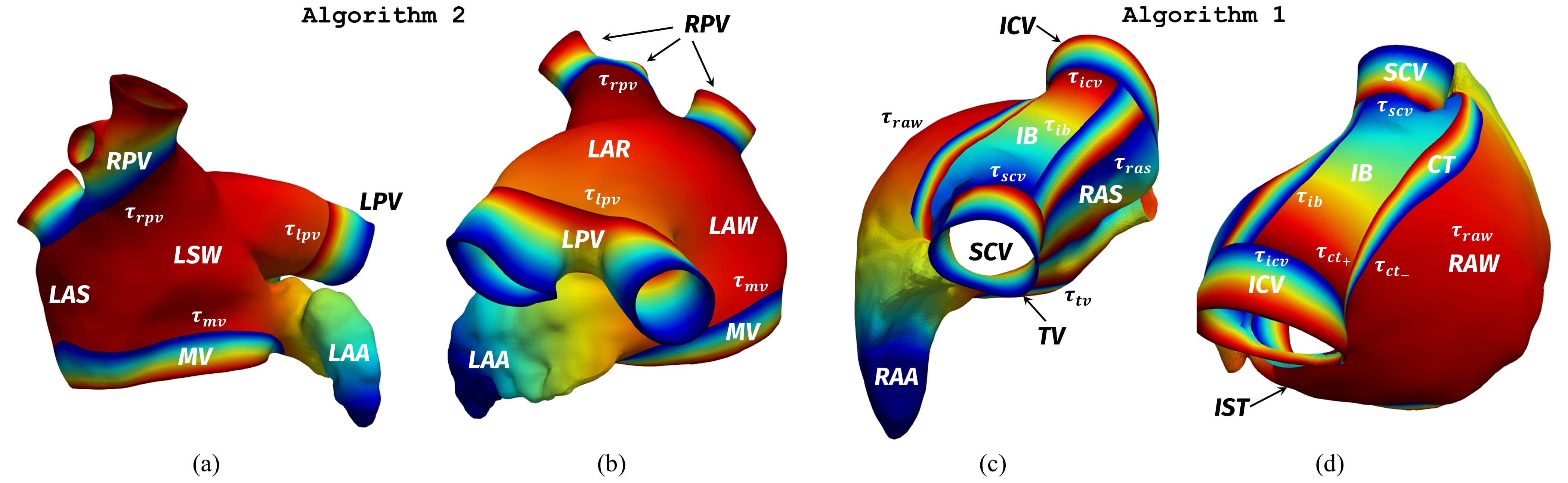}
	\caption{\textcolor{black}{Definition of the bundles and their dimension throughout the atrial domain following the rules reported in Algorithms~1 and 2; (a,b): bundles selection in LA following Algorithms 2; (c,d): bundles selection in RA following Algorithms~1. SCV, ICV: superior and inferior caval veins; LPV, RPV: left and right pulmonary veins; TV, MV: tricuspid and mitral valve rings; RAA, LAA: right and left appendage; RAS, LAS: right and left atrial septum; RAW, LAW: right and left atrial lateral wall; LSW: left septum wall; LAR: left atrial roof; IB: inter-caval bundle; CT: crista terminalis; IST: isthmus. $\tau_{tv}$, $\tau_{raw}$, $\tau_{ct_-}$, $\tau_{ct_+}$, $\tau_{icv}$, $\tau_{scv}$, $\tau_{ib}$ and $\tau_{ras}$ are the parameters related to the dimension of RA bundles, while $\tau_{mv}$, $\tau_{lpv}$ and $\tau_{rpv}$ to that of LA bundles.}}
	\label{fig:Algo12}
\end{figure}
\begin{algorithm}[ht!]
	\caption{{\sl Bundles selection for RA in the atrial LDRBM}}
	\noindent Let \textcolor{black}{$\tau_{tv}$, $\tau_{raw}$}, $\tau_{ct_-}$, $\tau_{ct_+}$, $\tau_{icv}$, $\tau_{scv}$, \textcolor{black}{$\tau_{ib}$ and $\tau_{ras}$} be the parameters referring to the \textcolor{black}{related} bundles.
	\begin{algorithmic}		
		\If {$\psi_r \ge$ \textcolor{black}{$\tau_{tv}$}} \textbf{set} $\boldsymbol{k}=\nabla \psi_r$ \textcolor{black}{$\longrightarrow$ TV}	
		\Else
		\If{$\psi_r <$ \textcolor{black}{$\tau_{raw}$}}
		\If{$\psi_w \ge \tau_{ct_-}$ and $\psi_w \le \tau_{ct_+}$}  \textbf{set} $\boldsymbol{k}=\nabla \psi_w$ \textcolor{black}{$\longrightarrow$ CT}
		\ElsIf{$\psi_w < \tau_{ct_-}$}		
		\If{$\psi_v \ge \tau_{icv}$ or $\psi_v \le \tau_{scv}$} \textbf{set} $\boldsymbol{k}=\nabla \psi_v$ \textcolor{black}{$\longrightarrow$ ICV and SCV}
		\Else $ \,$ \textbf{set} $\boldsymbol{k}=\nabla \psi_{ab}$ \textcolor{black}{$\longrightarrow$ RAW}
		\EndIf		
		\Else
		\If{$\psi_v \ge \tau_{icv}$ or $\psi_v \le \tau_{scv}$} \textbf{set} $\boldsymbol{k}=\nabla \psi_v$ \textcolor{black}{$\longrightarrow$ ICV and SCV}	 
		\Else
		\If{$\psi_w \le$ \textcolor{black}{$\tau_{ib}$}} \textbf{set} $\boldsymbol{k}=\nabla \psi_{v}$ \textcolor{black}{$\longrightarrow$ IB}
		\ElsIf{$\psi_w \ge$ \textcolor{black}{$\tau_{ras}$}} \textbf{set} $\boldsymbol{k}=\nabla \psi_{r}$ \textcolor{black}{$\longrightarrow$ RAS (centre)}
		\Else $\,$ \textbf{set} $\boldsymbol{k}=\nabla \psi_{w}$ \textcolor{black}{$\longrightarrow$ RAS (top)}
		\EndIf
		\EndIf
		\EndIf
		\Else
		\If{$\psi_v \ge \tau_{icv}$ or $\psi_v \le \tau_{scv}$} \textbf{set} $\boldsymbol{k}=\nabla \psi_v$ \textcolor{black}{$\longrightarrow$ ICV and SCV}	
		\Else 
		\If{$\psi_w \ge 0$} \textbf{set} $\boldsymbol{k}=\nabla \psi_{r}$ \textcolor{black}{$\longrightarrow$ RAS (bottom)}
		\Else $\,$ \textbf{set} $\boldsymbol{k}=\nabla \psi_{ab}$ \textcolor{black}{$\longrightarrow$ RAW, IST and RAA}
		\EndIf
		\EndIf
		\EndIf			
		\EndIf				
	\end{algorithmic}
	\label{RA-LDRBM}	
\end{algorithm}
\begin{algorithm}[t!]
	\caption{{\sl Bundles selection for LA in the atrial LDRBM}}
	\noindent Let \textcolor{black}{$\tau_{mv}$}, $\tau_{lpv}$ and $\tau_{rpv}$ be the parameters referring to the \textcolor{black}{related} bundles.
	\begin{algorithmic}
		\If{$\psi_r \ge$ \textcolor{black}{$\tau_{mv}$}} \textbf{set} $\boldsymbol{k}=\nabla \psi_r$ \textcolor{black}{$\longrightarrow$ MV}
		\Else
		\If {$\psi_v \ge \tau_{lpv}$ or $\psi_v \le \tau_{rpv}$} \textbf{set} $\boldsymbol{k}=\nabla \psi_v$ \textcolor{black}{$\longrightarrow$ LPV and RPV}
		\Else $\,$ \textbf{set} $\boldsymbol{k}=\nabla \psi_{ab}$ \textcolor{black}{$\longrightarrow$ LAS, LSW, LAW, LAR and LAA} 
		\EndIf
		\EndIf
	\end{algorithmic}
	\label{LA-LDRBM}
\end{algorithm}
\paragraph{3. Bundles selection}
Definition of the atrial bundles and their dimension throughout the domain $\Omega_{atrial}$, in order to match histology and DTI observations. With this aim, the atrial LDRBM assigns, for each point in $\Omega_{atrial}$, a unique \textcolor{black}{intra-atrial} distance $\psi_i$, among those defined in step 2, following the rules reported in Algorithms 1 and 2 for the bundle selection in the right and left atrium, respectively (see step 3 in Figure~\ref{fig:LARA-LDRBM} \textcolor{black}{and also Figure \ref{fig:Algo12}}). During this assignment, the atrial LDRBM defines a unique normal direction $\boldsymbol k$ by taking the gradient of a specific \textcolor{black}{intra-atrial} distances, $\boldsymbol k= \nabla \psi_i$. \textcolor{black}{Following Algorithms \ref{RA-LDRBM} and \ref{LA-LDRBM} the principal anatomical atrial regions are introduced: for RA, superior (SCV) and inferior caval veins (ICV), tricuspid valve ring (TV), right appendage (RAA), septum (RAS), inter-caval bundle (IB), crista terminalis (CT), isthmus (IST) and atrial later wall (RAW); for LA, left (LPV) and right pulmonary veins (RPV), mitral valve ring (MV), left appendage (LAA), septum (LAS), septum wall (LSW), atrial lateral wall (LAW) and roof (LAR), see Figure \ref{fig:Algo12}.} \textcolor{black}{Moreover, in order to specify the bundles dimension, the parameters $\tau_{i}$ are introduced: for the RA $\tau_{tv}$, $\tau_{ib}$, $\tau_{icv}$, $\tau_{scv}$, $\tau_{ct_+}$, $\tau_{ct_-}$, $\tau_{raw}$ and $\tau_{ras}$ refer to  
TV, IB, ICV, SCV, and upper and lower limit of CT  
bundles, respectively; 
for LA $\tau_{mv}$, $\tau_{lpv}$ and $\tau_{rpv}$ refer to 
MV, LPV and RPV 
bundles, respectively (see Figure \ref{fig:Algo12}).} The complete bundles selection procedures are fully displayed in Algorithms \ref{RA-LDRBM} and \ref{LA-LDRBM} (see also step 3 in Figure \ref{fig:LARA-LDRBM} \textcolor{black}{and Figure~\ref{fig:Algo12}}).   

\paragraph{4. Local coordinate system}
Definition of the myofiber orientations by means of an orthonormal local coordinate system, built at each point of the atrial domain. This step is performed in the same way as for the ventricles: the gradient of the transmural distance $\phi$ is used to build the transmural direction $\nabla \phi$ which is taken as one input of the function \texttt{axis} \eqref{axis} together with the unique normal direction $\boldsymbol k$:  
\begin{equation}
\label{coo_atrial}
Q=[\widehat{\boldsymbol{e}}_l, \widehat{\boldsymbol{e}}_n, \widehat{\boldsymbol{e}}_t]
=[\boldsymbol f, \boldsymbol n, \boldsymbol s]=\texttt{axis}(\boldsymbol k, \nabla \phi),
\end{equation}
where $\widehat{\boldsymbol{e}}_l$, $\widehat{\boldsymbol{e}}_n$ and $\widehat{\boldsymbol{e}}_t$ are the unit longitudinal, normal and transmural directions, respectively. Moreover, since we are not prescribing any transmural variation in the fiber bundles, the three unit directions correspond to the final fiber, sheet and cross-fiber directions $\boldsymbol f$, $\boldsymbol n$ and $\boldsymbol s$ (see step 4 in Figure \ref{fig:LARA-LDRBM}). $\blacksquare$

\paragraph{\textcolor{black}{Atrial LDRBM rules}}
\textcolor{black}{Algorithms 1 and 2, combined with the definition of the local coordinate system \eqref{coo_atrial}, prescribe a fiber field $\boldsymbol f$ based on the following rules derived from histo-anatomical and DT-MRI fiber data observations \cite{pashakhanloo2016myofiber,sanchez2014left,ho2002atrial,ho2009importance,ho2012left,sanchez2013standardized,roney2019universal}, see  Figures~\ref{fig:Algo12}:}

\begin{enumerate}
	\item[\textcolor{black}{R1:}] \textcolor{black}{Circular fiber arrangements are prescribed around LPV, RPV, SCV, ICV, TV, MV, and encircle both appendages (RAA and LAA) \cite{pashakhanloo2016myofiber,sanchez2014left,roney2019universal};} 
	\item[\textcolor{black}{R2:}] \textcolor{black}{Fibers direction of CT runs longitudinally from the base of SCV to ICV \cite{sanchez2015anatomical};}
	\item[\textcolor{black}{R3:}] \textcolor{black}{RA structures like IB and RAW are almost vertically oriented, whereas those of RAS are parallel to CT \cite{ho2002atrial,ho2009importance};}
	\item[\textcolor{black}{R4:}] \textcolor{black}{IST fibers have the same direction of those of TV \cite{sanchez2013standardized};}
	\item[\textcolor{black}{R5:}] \textcolor{black}{LAS fibers are aligned (parallel) to the nearby RAS  \cite{ferrer2015detailed};} 
	\item[\textcolor{black}{R6:}] \textcolor{black}{Directions of LAR and LAW descend perpendicularly to MV, while fibers of LSW present a smooth transition going to LAS and LAA \cite{ferrer2015detailed,ho2012left,roney2019universal}.}
\end{enumerate}   

\section{Modelling cardiac electrophysiology}\label{sec:num}
\textcolor{black}{In this section we briefly recall the mathematical model for the description of the EP activity in the cardiac tissue, that is the monodomain equation endowed with suitable ionic models for human action potential \cite{quarteroni2017integrated,dossel2012computational,franzone2014mathematical,vigmond2003computational,plank2008mitochondrial,trayanova2014exploring,niederer2019computational}.}

Cardiac tissue is an orthotropic material, arising from the cellular organization of the myocardium in fibers, laminar sheet and cross-fibers, which is mathematically modelled by the conductivity tensor    
\begin{equation}
\label{D_tensor}
\mathbb{D} = \sigma_f \boldsymbol f \otimes \boldsymbol f + \sigma_s \boldsymbol s \otimes \boldsymbol s + \sigma_n \boldsymbol n \otimes \boldsymbol n,
\end{equation}
where $\sigma_f$, $\sigma_s$ and $\sigma_n$ are the conductivities along fibers ($\boldsymbol f$), sheets ($\boldsymbol s$), and cross-fibers ($\boldsymbol n$) directions, respectively.
Given a computational domain \(\Omega\) and a time interval \((0, T]\), the monodomain equations read:\\
find, for each $t$, the transmembrane action potential $u:\Omega \times (0, T]\to\mathbb R$ and the gating variables $\boldsymbol w:\Omega \times (0, T]\to\mathbb{R}^n$, such that
\begin{subequations}
	\label{eq:monodomain_system}
	\begin{align}
	& \chi C_m\frac{\partial u}{\partial t} - \nabla\cdot\left(\mathbb{D}
	\nabla u\right) + \chi I_{ion}(u, \boldsymbol w) = I_{app}(\mathbf{x},t) & 
	\text{in } \Omega \times (0,T],  \label{eq:monodomain}\\ 
	& \displaystyle\frac{d\boldsymbol w}{dt} = \boldsymbol G(u,\boldsymbol w) \label{eq:ionicmodel} & \text{in } \Omega \times (0,T] ,
	\end{align}
\end{subequations} 
where $\chi$ is the surface area-to-volume ratio of cardiomyocytes, $C_m$ is the specific trans-membrane capacitance per unit area, $I_{app}$ is an external applied  current which serves to initiate the signal propagation, $I_{ion}$ and $\boldsymbol G \in \mathbb{R}^n$ are the reaction terms, linking the macroscopic action potential propagation to the cellular dynamics. 
The unknown $\boldsymbol w$ is a $n$--th dimensional vector function fulfilling a system of differential algebraic equations representing the percentage of open channels per unit area of the membrane.
Specifically, we used the Courtemanche-Ramirez-Nattel (CRN, $n=20$) in case of atrial action potential and the Ten-Tusscher-Panfilov (TTP, $n=18$) for the ventricular one (for further details see \cite{courtemanche1998ionic} for CRN and \cite{ten2006alternans} for TTP). 
Furthermore, system \eqref{eq:monodomain_system} is equipped with suitable initial conditions for $u$ and $\boldsymbol w$ and homogeneous Neumann boundary conditions for $u$ at the boundary $\partial\Omega$.

\textcolor{black}{For the time discretization of the monodomain system \eqref{eq:monodomain_system} we consider a {\sl Backward Difference Formulae} approximation of order $3$ (BDF3) with an explicit treatment of the reaction term. Moreover, the diffusion term is treated implicitly, whereas the ionic terms explicitly \cite{franzone2014mathematical,plank2008mitochondrial,vigmond2008solvers}.}

\textcolor{black}{Regarding the space discretization, we used continuous Finite Elements (FE) on hexahedral meshes and the discretization of the ionic current term \(I_{ion}\) is performed following the Ionic Current Interpolation approach \cite{trayanova2011whole,quarteroni2017integrated,hurtado2018non,jilberto2018semi,arevalo2016arrhythmia,krishnamoorthi2013numerical}.}
	


\section{Numerical results}\label{sec:res}
This section is dedicated to several results both for the fibers generation and the numerical EP simulations. These have been performed both on idealized and realistic human ventricular and atrial models. As realistic geometry, we use the 3D heart model {\sl Zygote} \cite{zygote2014}, a complete heart geometry reconstructed from high-resolution CT-scans representing an average healthy heart. Being a very detailed model of the human heart, it demonstrates the applicability of the proposed methods to arbitrary patient-specific scenarios.

We organize the section as follows. After a brief description related to the setting of numerical simulations (Section \ref{sec:setting}),
we show various comparisons among the three LDRBMs for ventricles fiber generation (Section~\ref{sec:ventricles}): we compare the fiber orientations and we analyse their influence in terms of activation times computed as output of numerical EP simulations\footnote{The activation time of a given point in the cardiac muscle is defined as the time when the transmembrane potential derivative $\frac{du}{dt}$ reaches its maximal value.}. For this comparison, first we make use of an idealized biventricular geometry built using the prolate spheroid coordinate systems \cite{sermesant2005simulation}, and then we consider the Zygote biventricular model \cite{zygote2014}. Section \ref{sec:atria} is devoted to the novel atrial LDRBM. We show fiber bundles reconstruction applied to an idealized case \cite{pegolotti2019isogeometric}, to the realistic Zygote geometries \cite{zygote2014} \textcolor{black}{and to the atrial geometry presented in e.g. \cite{tobon2013three,ferrer2015detailed}, in what follows referred to as {\sl Riunet geometry} \textcolor{black}{\footnote{\textcolor{black}{Freely available online at \url{https://riunet.upv.es/handle/10251/55150}.}}} (from the name of the repository).}  
Afterwards, we investigate the influence of atrial fibers in EP simulations comparing the fiber activation map with respect to a one obtained with an isotropic conductivity. Finally, in Section \ref{sec:heart} we present an EP simulation of a realistic four chamber heart including fibers generated by LDRBMs for both atria and ventricles.

\subsection{Setting of numerical simulations}\label{sec:setting}

\textcolor{black}{In order to build FE meshes, a pre-processing phase was applied to every ventricular and atrial geometry used in the simulations. For this preprocessing phase we deeply rely on the novel semi-automatic meshing tool proposed in \cite{fedele2019processing} which consists of multiple steps including tagging, geometry smoothing and hexahedral FE mesh generation. Specifically, the Vascular Modelling Toolkit \texttt{vmtk} software \cite{antiga2008vascular} (\url{http://www.vmtk.org}) together with a new meshing tool were used to perform such pre-processing phase. The tagging procedures carried out in this work, for the ventricular and atrial LDRBMs, are detailed in the Appendix. We remark however that our tagging procedure is not crucial for the applicability of the ventricular and atrial LDRBMs. Indeed, they are perfectly compatible with other tagging processes presented in other works (see for example \cite{doste2019rule} for the ventricles and \cite{strocchi2020simulating} for the whole heart case).}

The numerical approximation of the monodomain system \eqref{eq:monodomain_system} requires the following physical data: the transmembrane capacitance per unit area $C_m$, the membrane surface-to-volume ratio $\chi$ and the conductivities along the three direction of the myofibers $\sigma_f$, $\sigma_s$ and $\sigma_n$ appearing in the conductivity tensor $\mathbb{D}$. The values chosen for the first two quantities are $C_m=1$ $\mu$F/$\mathrm{cm}^{2}$ and $\chi=1400$ $\mathrm{cm}^{-1}$
\cite{potse2006comparison,quarteroni2019,roth1991action,niederer2011verification}. The conductivity values $\sigma_f$, $\sigma_s$ and $\sigma_n$ were fitted by an iterative procedure \textcolor{black}{(reported in the Appendix)} fully described in \cite{costa2013automatic} (see also \cite{fastl2018personalized,augustin2016anatomically}) in order to match the following conduction velocity values: for the ventricles, $60$ cm/s in the fiber direction $\boldsymbol f$,  $40$ cm/s in the sheet direction $\boldsymbol s$ and $20$ cm/s in the normal direction $\boldsymbol n$ \cite{augustin2016anatomically}; for the atria, $120$ cm/s in the fiber direction $\boldsymbol f$ and $40$ cm/s along the sheet $\boldsymbol s$ and cross-fiber directions $\boldsymbol n$ \cite{fastl2018personalized, augustin2019impact, monaci2018computational}. Finally, to initiate the signal propagation in the cardiac muscle, the monodomain system \eqref{eq:monodomain_system} requires to specify the external applied current $I_{app}(\mathbf{x},t)$. In this work $I_{app}(\mathbf{x},t)$ was modeled as a series of impulses (with radius $2.5$ mm and duration $3$ ms) applied in spherical subsets of the domain and prescribed alongside the \textcolor{black}{ventricular} and \textcolor{black}{atrial endocardia}. Its amplitude is $50000$ $\mu$A/$\mathrm{cm}^{3}$, for both \textcolor{black}{atrial} and \textcolor{black}{ventricular} domains, in agreement with \cite{niederer2011verification}.
We used this value for all the simulations, while the stimuli locations will be specified for each case reported in Sections \ref{sec:ventricles},  \ref{sec:atria} and \ref{sec:heart}.

Regarding the mesh element size $h$ and the time step $\Delta t$, related to the
space and time discretizations of the system \eqref{eq:monodomain_system}
we used continuous FE of order 1 ($Q_1$) on hexahedral meshes with an average mesh size of $h=350$ $\mu$m
\cite{trayanova2011whole,quarteroni2017integrated,plank2008mitochondrial,vigmond2008solvers,hurtado2018non,jilberto2018semi,arevalo2016arrhythmia,niederer2011verification}. 
and $\text{BDF}$ of order $\sigma=3$ with a time step of $\Delta t=50$ $\mu$s \textcolor{black}{(see the Appendix for further details).} 

All the \textcolor{black}{ventricular} LDRBMs (presented in Section \ref{sec:ventr}), the novel atrial LDRBM (detailed in Section \ref{sec:atrial}) and the numerical methods for the EP presented in Section \ref{sec:num} have been implemented within \texttt{life\textsuperscript{x}} (\url{https://lifex.gitlab.io/lifex}), a new in-house high-performance \texttt{C++} FE library mainly focused on cardiac applications based on \texttt{deal.II} FE core \cite{dealII91} (\url{https://www.dealii.org}). All the numerical simulations were executed on the cluster \texttt{iHeart} (Lenovo SR950 
8x24-Core Intel Xeon Platinum 8160, 2100 MHz and 1.7TB RAM) at MOX, 
Dipartimento di Matematica, Politecnico di Milano.   

To analyse the results we used \texttt{ParaView} (\url{https://www.paraview.org}) an open-source, multi-platform data analysis and visualization application. In particular, to visualize the fiber fields we applied in sequence the \texttt{streamtracer} and the \texttt{tube} \texttt{ParaView} filters.


\subsection{Ventricular fibers generation and electrophysiology}\label{sec:ventricles}
This section is dedicated to the comparison among the different fiber fields generated by the three ventricular LDRBMs presented in Section \ref{sec:ventr}. Moreover, we also investigate the influence of the different fibers orientation in the activation times, produced by EP simulations. We perform both the analyses first in an idealized biventricular geometry \cite{razumov2018study} and then in the Zygote realistic human model \cite{zygote2014}.

\subsubsection{Idealized biventricular model}
The first comparison among the three \textcolor{black}{ventricular} LDRBMs was performed on a well established idealized biventricular geometry that has been used in several computational studies \cite{wong2014generating,sermesant2005simulation,razumov2018study,goktepe2010atrial,goktepe2010electromechanics,ahmad2018multiphysics}
and for \textcolor{black}{ventricular} volume estimation from 2D images \cite{mercier1982two}. The heart ventricles are approximated as two intersecting truncated ellipsoids. 

We constructed the idealized biventricular geometry using the prolate spheroid coordinate systems in the built-in CAD engine of \texttt{gmsh}, an open source 3D finite element mesh generator (\url{http://gmsh.info}), see Figure \ref{fig:BiV_Ideal_Fibers}. For the details about the geometrical definition of the idealized biventricular, we refer to~\cite{razumov2018study,ahmad2018multiphysics}. 

\subsubsection*{5.2.1.1 Fiber generation in the idealized ventricles}
Fiber orientations obtained for the three LDRBMs (R-RBM, B-RBM and D-RBM) in the idealized biventricular model are shown in Figures \ref{fig:BiV_Ideal_Fibers}(a-f). 
The input angles values $\alpha_{endo,\ell}$, $\alpha_{epi,\ell}$, $\alpha_{endo,r}$, $\alpha_{epi,r}$, $\beta_{endo,\ell}$, $\beta_{epi,\ell}$, $\beta_{endo,r}$ and $\beta_{epi,r}$ were chosen for all the three methods based on the observations of histological studies in the human heart \cite{lombaert2012human,greenbaum1981left,anderson2009three,ho2006anatomy,sanchez2015anatomical,lunkenheimer2013models,stephenson2016functional} (see also \cite{doste2019rule}):
\begin{figure}[t!]	
	\centering
	\includegraphics[width=0.9\textwidth]{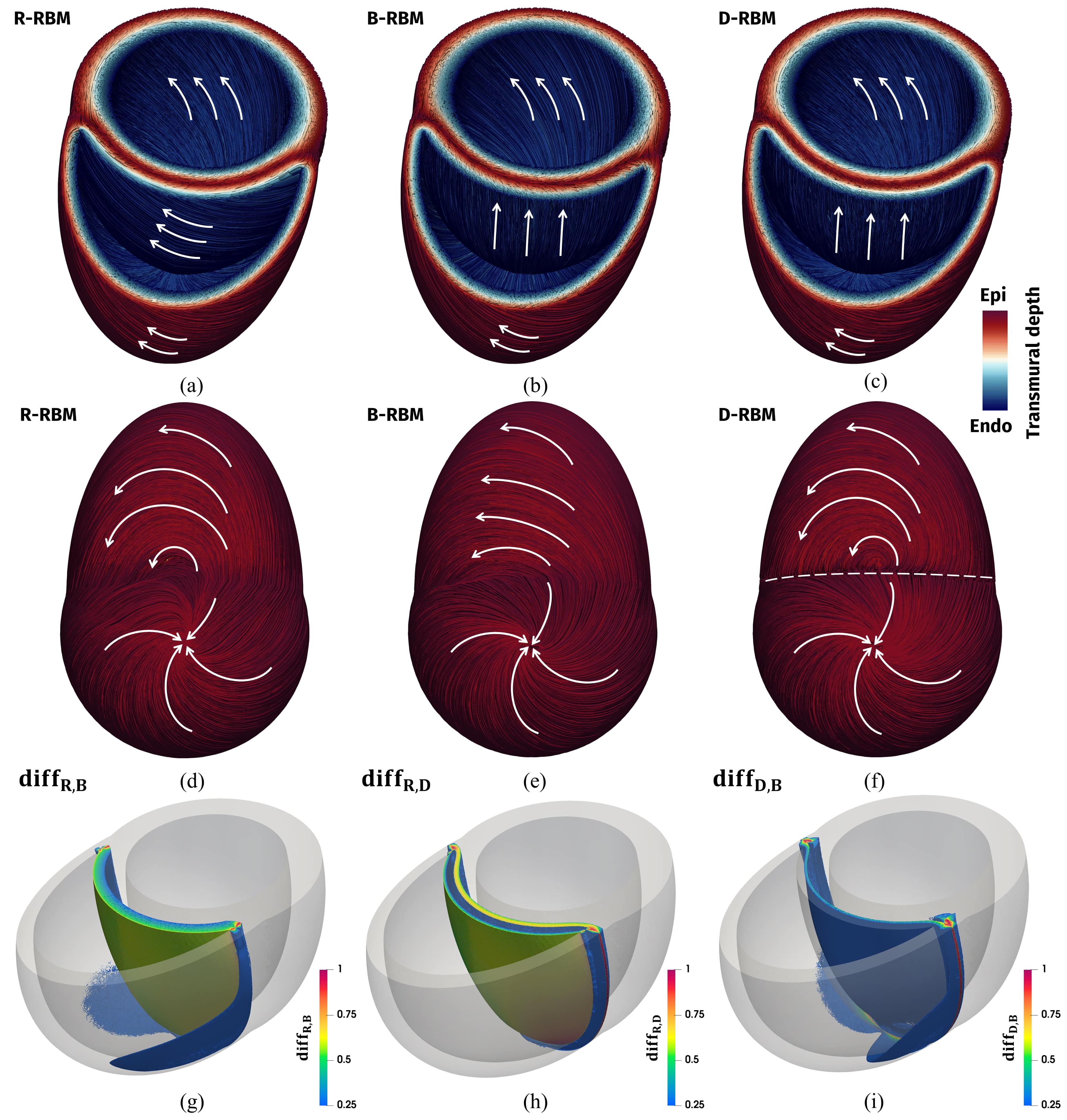}
	\caption{Comparison among LDRBMs in the idealized biventricular model. Streamlines of the vector $\boldsymbol f$ is depicted for R-RBM (a,d), B-RBM (b,e) and D-RBM (c,f). \textcolor{black}{White arrows represent the main fibers direction in specific ventricular regions (displayed in red for the epicardium and in blue for the endocardium); the dashed line in Figure (f) highlights the inter-ventricular junctions discontinuity of D-RBM}. Top: Frontal view; Centre: apex view; Bottom: Differences $\mathrm{diff}_{i,j}$ among the three LDRBMs, $\mathrm{diff}_{\mathrm{R,B}}$ (g), $\mathrm{diff}_{\mathrm{R,D}}$ (h) and  $\mathrm{diff}_{\mathrm{D,B}}$ (i); only values $\mathrm{diff}_{i,j} \ge 0.25$ are displayed.}
	\label{fig:BiV_Ideal_Fibers}
\end{figure}
\begin{equation}
\label{angles}
\begin{alignedat}{2}
& \alpha_{epi,\ell}=-60^o, \quad \alpha_{endo,\ell}=+60^o, \quad \alpha_{epi,r}=-25^o, \quad \alpha_{endo,r}=+90^o;
\\
& \beta_{epi,\ell}=+20^o, \quad \beta_{endo,\ell}=-20^o, \quad \beta_{epi,r}=+20^o, \quad \beta_{endo,r}=0^o.
\end{alignedat}
\end{equation}  

We observe that all the LDRBMs represent the characteristic helical structure of the left ventricle and a compatible fiber orientations both in the right endocardium, not facing to the septum, and in the right epicardium, far enough from the inter-ventricular junctions. Most of the differences occur in the right \textcolor{black}{ventricular} endocardium facing the septum (see Figures \ref{fig:BiV_Ideal_Fibers}(a-c)), in the inter-ventricular junctions between the two ventricles \textcolor{black}{and in the right epicardial lower region} (see Figures \ref{fig:BiV_Ideal_Fibers}(d-f)). 

We computed the difference $\mathrm{diff}_{i,j}$ of the fiber field $\boldsymbol{f}$ among the three methods, defined as:  
\begin{equation}
\label{diff_F}
\mathrm{diff}_{i,j}(\boldsymbol{x})=1-|\boldsymbol{f}_i(\boldsymbol{x}) \cdot \boldsymbol{f}_j(\boldsymbol{x})| \qquad \qquad i,j=\mathrm{R,B,D} \quad (i \ne j), 
\end{equation} 
where $\boldsymbol{f}_{\mathrm{R}}$, $\boldsymbol{f}_{\mathrm{B}}$ and $\boldsymbol{f}_{\mathrm{D}}$ are the vector fiber fields of R-RBM, B-RBM and D-RBM, respectively. If $\boldsymbol{f}_i$ and $\boldsymbol{f}_j$ are parallel, $\mathrm{diff}_{i,j}=0$, otherwise, when orthogonal, $\mathrm{diff}_{i,j}=1$. The result of these comparisons is reported in Figures \ref{fig:BiV_Ideal_Fibers}(g), \ref{fig:BiV_Ideal_Fibers}(h) and \ref{fig:BiV_Ideal_Fibers}(i). As expected, most of the discrepancies are in the septum, in the inter-ventricular junctions \textcolor{black}{and in the right epicardial lower region.}
\begin{figure}[t!]	
	\centering
	\includegraphics[width=0.9\textwidth]{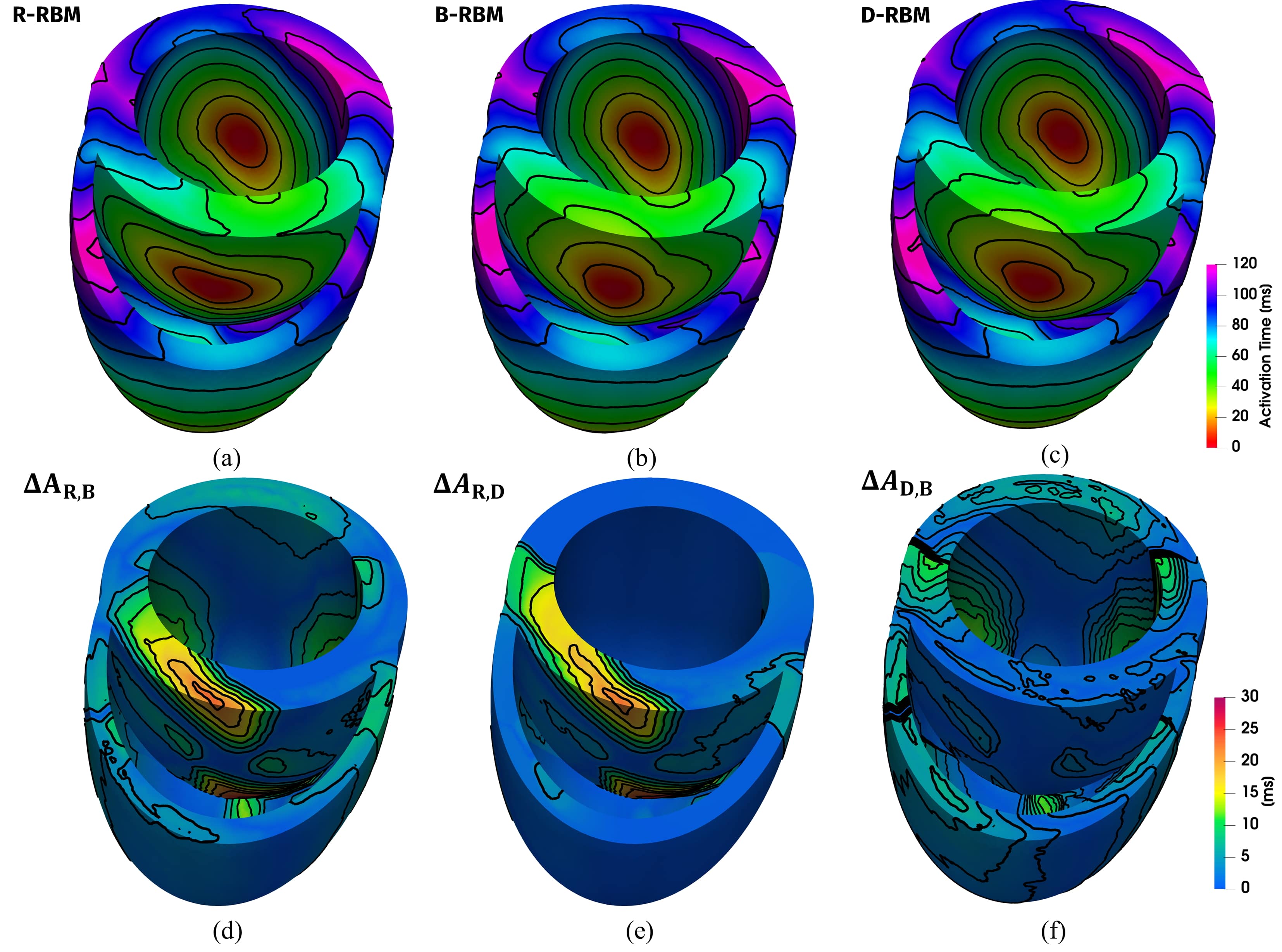}
	\caption{Top: Activation time for R-RBM (a), B-RBM (b) and D-RBM (c) in the idealized biventricular model. Bottom: Absolute difference among the activation maps, $\Delta A_{\mathrm{R,B}}$ (d), $\Delta A_{\mathrm{R,D}}$ (e) and $\Delta A_{\mathrm{D,B}}$ (f).}
	\label{fig:BiV_Ideal_AT}
\end{figure}
\subsubsection*{5.2.1.2 Electrophysiology in the idealized ventricles}
In order to evaluate the influence of the three LDRBMs fiber architectures in the electric signal propagation through the cardiac muscle, we performed three EP simulations (with the setting detailed in Section~\ref{sec:param}), one for each LDRBM. To initiate the action potential propagation we applied four endocardial stimuli: two for each ventricle, one in the mid-septal zone and one in the \textcolor{black}{lateral endocardial wall.} In Figures \ref{fig:BiV_Ideal_AT}(a-c) we report the activation maps obtained with the three fibers configurations. The activation pattern for all the three methods are very similar in the left and right ventricles, while most of the differences are visible in the septum, see Figures \ref{fig:BiV_Ideal_AT}(a), \ref{fig:BiV_Ideal_AT}(b) and \ref{fig:BiV_Ideal_AT}(c). We computed also the absolute difference $\Delta A_{i,j}(\boldsymbol{x})$ in the activation pattern among the different methods as:  
\begin{equation}
\label{difference}
\Delta A_{i,j}(\boldsymbol{x})=|A_{i}(\boldsymbol{x})-A_j(\boldsymbol{x})| \qquad \qquad i,j=\mathrm{R,B,D} \quad (i \ne j), 
\end{equation}  
where $A_{\mathrm{R}}$, $A_{\mathrm{B}}$ and $A_{\mathrm{D}}$ are the activation times for R-RBM, B-RBM and D-RBM, respectively (see Figures~\ref{fig:BiV_Ideal_AT}(d), \ref{fig:BiV_Ideal_AT}(e) and \ref{fig:BiV_Ideal_AT}(f)). 

The most remarkable differences in both $\Delta A_{\mathrm{R,B}}$ and $\Delta A_{\mathrm{R,D}}$ are exhibited in the septum, particularly in the part facing the right endocardium, while $\Delta A_{\mathrm{B,D}}$ never exceeds $15$ ms, see Figures~\ref{fig:BiV_Ideal_AT}(d), \ref{fig:BiV_Ideal_AT}(e) and \ref{fig:BiV_Ideal_AT}(f). Also in the activation maps, as \textcolor{black}{expected}, we retrieve differences in the septum zone caused by the different fiber orientations on that region made by the three methods, as seen in the fibers comparison, see \textcolor{black}{Figure} \ref{fig:BiV_Ideal_Fibers}.

Finally, we evaluated the maximal discrepancies, $M_{i,j}=\max_{\boldsymbol{x}\in \Omega_{myo}}\Delta A_{i,j}(\boldsymbol{x})$, among the three methods, which are: 
\begin{equation*}
M_{\mathrm{R,B}}=35\,\mathrm{ms}, \qquad M_{\mathrm{R,D}}=33\,\mathrm{ ms}, \qquad M_{\mathrm{B,D}}=15\,\mathrm{ ms}.  
\end{equation*}  
The location of both $M_{\mathrm{R,B}}$ and $M_{\mathrm{R,D}}$ is in the lower part of the right \textcolor{black}{ventricular} septum, while $M_{\mathrm{B,D}}$ is placed in the lower anterior region of the left ventricle, see Figures \ref{fig:BiV_Ideal_AT}(d), \ref{fig:BiV_Ideal_AT}(e) and \ref{fig:BiV_Ideal_AT}(f). Considering a total activation time of about $A_{max}=120$ ms for the all biventricular muscle, the maximum relative differences, $M^{\%}_{i,j}=M_{i,j}/A_{max}$, among the three \textcolor{black}{ventricular} LDRBMs are  
\begin{equation*}
M^{\%}_{\mathrm{R,B}}=29 \%, \qquad M^{\%}_{\mathrm{R,D}}=28 \%, \qquad M^{\%}_{\mathrm{B,D}}=13 \%.  
\end{equation*} 

\subsubsection{Realistic full biventricular model}\label{sec:real-biv}
\begin{figure}[t!]	
	\centering
	\includegraphics[width=0.825\textwidth]{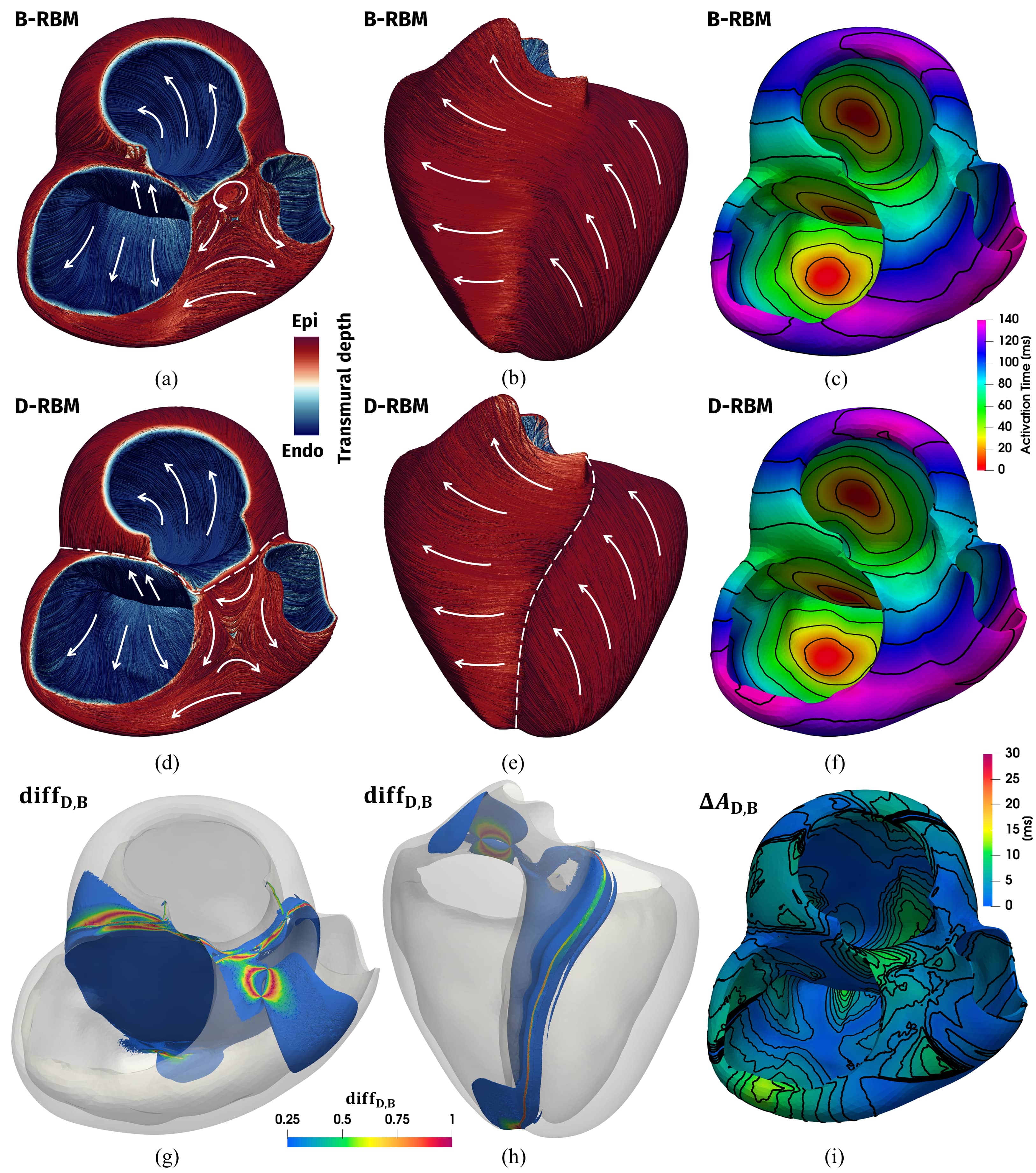}
	\caption{Comparison for B-RBM and D-RBM in a realistic full biventricular model. Top (a-c): B-RBM. Centre (d-f): D-RBM. Bottom (g-i): differences between B-RBM and D-RBM. Streamlines of the vector $\boldsymbol f$: top (a,d) and lateral views (b,e). \textcolor{black}{White arrows represent the main fibers direction in specific ventricular regions (displayed in red for the epicardium and in blue for the endocardium); dashed line in Figures (d,e) highlights the inter-ventricular junctions discontinuity of D-RBM}. Difference in the fiber orientations $\mathrm{diff}_{\mathrm{D,B}}$ (g,h), only the values $\mathrm{diff}_{i,j} \ge 0.25$ are displayed. Activation maps using B-RBM and D-RBM: B-RBM (c) and D-RBM (f). Absolute difference among B-RBM and D-RBM activation maps, $\Delta A_{\mathrm{D,B}}$ (i).}
	\label{fig:BiV_Zygote_Fibers}
\end{figure}

The second comparison among the ventricular LDRBMs was performed on a realistic full biventricular geometry. For this purpose, we used the Zygote solid 3D heart model \cite{zygote2014}. In order to obtain a smooth endocardium in both ventricles, we removed all the papillary muscles and the trabeculae carneae, using the CAD modeller \texttt{SolidWorks} (\url{https://www.solidworks.com}) in combination with the software \texttt{Meshmixer}  (\url{http://www.meshmixer.com}), see Figure \ref{fig:BiV_Zygote_Fibers}. Considering the characteristics of the electrical signal propagation, and the anatomical constituents of the valvular and sub-valvular apparatus, we expect our calculations should not be substantially influenced by the papillary muscles elimination.

According to the motivations highlighted at the end of Section \ref{sec:ventr}, we performed a comparison in the full biventricular model only between B-RBM and D-RBM. 

\subsubsection*{5.2.2.1 Fibers generation in the realistic ventricles}
Fiber orientation for B-RBM and D-RBM in the Zygote full biventricular model are displayed in Figures~\ref{fig:BiV_Zygote_Fibers}(a-b) and \ref{fig:BiV_Zygote_Fibers}(d-e). We prescribed the same input angle values used for the ideal geometry, reported in \textcolor{black}{equation}~\eqref{angles}. Moreover, for D-RBM we have also specified the angles in the OT regions as follows \cite{doste2019rule}:
\begin{equation}
\label{OT-angles}
\alpha_{epi,OT}=0^o, \quad \alpha_{endo,OT}=+90^o, \quad \beta_{epi,OT}=0^o, \quad \beta_{endo,OT}=0^o.
\end{equation}  

The two LDRBMs well reproduce the helical structure of the left ventricle up to the mitral valve ring and exhibit a similar fiber orientation pattern in whole cardiac muscle, apart from the region between the tricuspid, the pulmonary and the aortic valve rings and far enough from the inter-ventricular junctions, see Figures \ref{fig:BiV_Zygote_Fibers}(a-b) and \ref{fig:BiV_Zygote_Fibers}(d-e). B-RBM presents a roll up in the fiber directions just after the aortic valve ring, while D-RBM has a more longitudinal fiber orientations in that region, see Figure \ref{fig:BiV_Zygote_Fibers}(a) and \ref{fig:BiV_Zygote_Fibers}(d). As also observed in the idealized case, the B-RBM fiber field in the inter-ventricular junctions has a smooth change passing from the left to the right ventricle, whereas D-RBM produces a strong discontinuity in the transition across the two ventricles, see Figure \ref{fig:BiV_Zygote_Fibers}(b) and \ref{fig:BiV_Zygote_Fibers}(e). 

We evaluated the mismatch of the fiber fields $\mathrm{diff}_{\mathrm{D,B}}$, defined in \eqref{diff_F}, between B-RBM and D-RBM.  
Indeed, $\mathrm{diff}_{\mathrm{D,B}}$ 
highlights the most relevant differences of the two methods in the septum, in the inter-ventricular junctions, in the regions of tricuspid, pulmonary and aortic valve rings and around the right \textcolor{black}{ventricular} apex, see Figures \ref{fig:BiV_Zygote_Fibers}(g) and \ref{fig:BiV_Zygote_Fibers}(h).

\subsubsection*{5.2.2.2 Electrophysiology in the realistic ventricles}
We performed two EP simulations (with the setting detailed in Section \ref{sec:param}), one with~B-RBM and one with D-RBM. Two stimuli were here applied to each ventricle: one in the mid-septal zone and one in the \textcolor{black}{lateral endocardial wall}. Figures \ref{fig:BiV_Zygote_Fibers}(c) and \ref{fig:BiV_Zygote_Fibers}(f) depict the computed activation times which result very similar in the whole myocardium. 
Figure \ref{fig:BiV_Zygote_Fibers}(i) shows the absolute difference between the two activation maps, $\Delta A_{\mathrm{D,B}}(\boldsymbol{x})=|A_{\mathrm{D}}(\boldsymbol{x})-A_{\mathrm{B}}(\boldsymbol{x})|$, where $A_{\mathrm{B}}$ and $A_{\mathrm{D}}$ are the activation times for B-RBM and D-RBM, respectively. We observe some discrepancies in the activation pattern near the tricuspid, the aortic and the pulmonary valve rings, and also in the endocardium near the right \textcolor{black}{ventricular} apex, although $\Delta A_{\mathrm{D,B}}$ never exceeds the value $14$ ms, see Figure \ref{fig:BiV_Zygote_Fibers}(i). The maximal relative discrepancy among the two methods is  $M_{\mathrm{D,B}}=\max_{\boldsymbol{x}\in \Omega_{myo}}\Delta A_{\mathrm{D,B}}(\boldsymbol{x})=14$ ms, corresponding to $M^{\%}_{\mathrm{D,B}}=M_{\mathrm{D,B}}/A_{max}=10 \, \%$,
with $A_{max}=140$ ms the total activation time.
The location of $M_{\mathrm{D,B}}$ is in the lower part of the endocardium just above the right \textcolor{black}{ventricular} apex, see Figure \ref{fig:BiV_Zygote_Fibers}(i). 

\begin{figure}[t!]	
	\centering
	\includegraphics[width=1\textwidth]{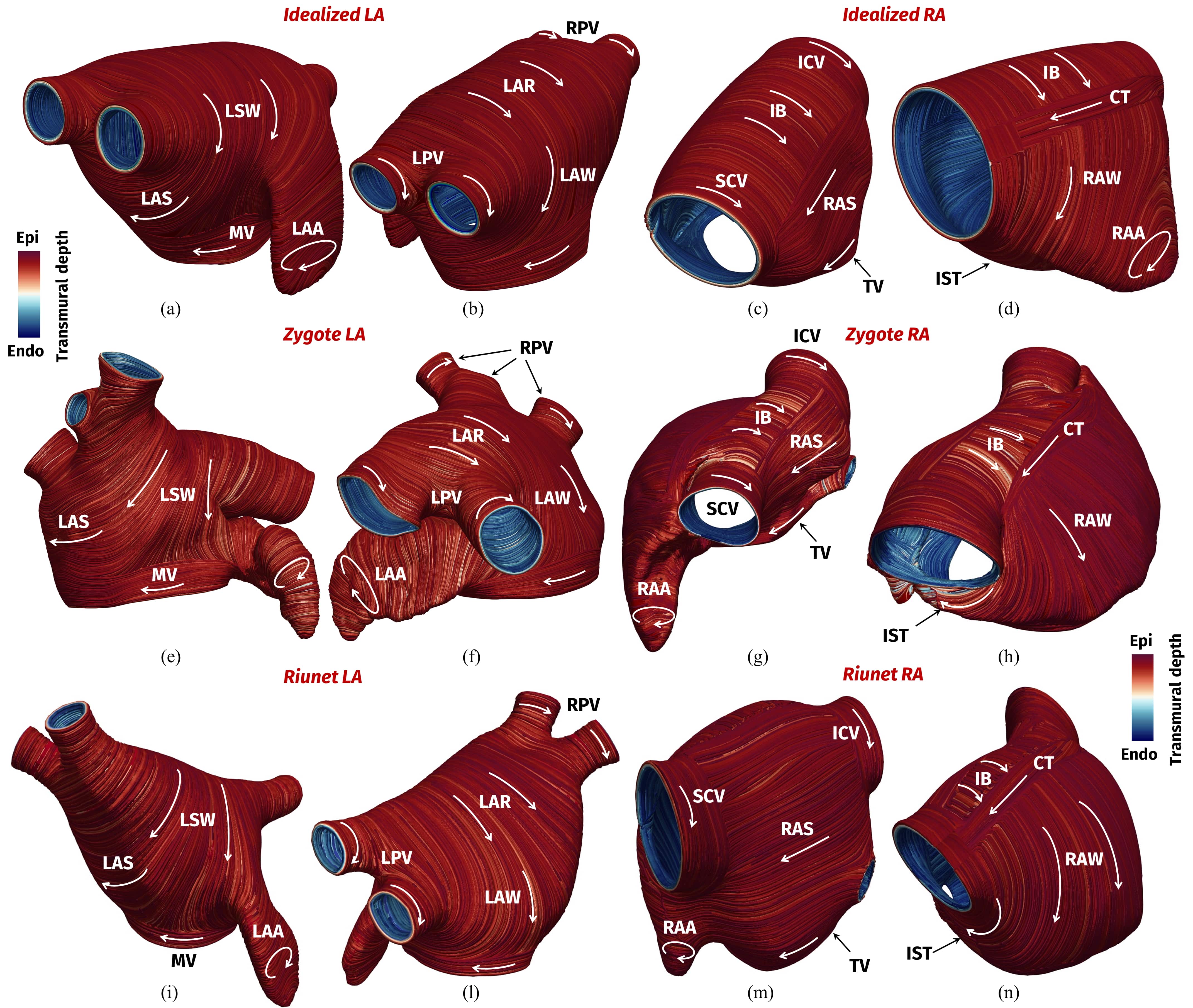}
	\caption{Atrial LDRBM fiber generation applied to idealized (a-d) and realistic (e-n) models \textcolor{black}{(Zygote (e-h) and Riunet (i-n) geometries).} Frontal (a,c,e,g,i,m) and dorsal (b,d,f,h,l,n) views of the atria. SCV, ICV: superior and inferior caval veins; LPV, RPV: left and right pulmonary veins; TV, MV: tricuspid and mitral valve rings; RAA, LAA: right and left appendage; RAS, LAS: right and left septum; RAW, LAW: right and left atrial lateral wall; LSW: left septum wall; LAR: left atrial roof; IB: inter-caval bundle; CT: crista terminalis; IST: isthmus.}
	\label{fig:Atria_Fibers}
\end{figure}  
\subsection{Atrial fibers generation and electrophysiology}\label{sec:atria}
We applied our new atrial LDRBM (presented in Section \ref{sec:atrial}) to reconstruct left and right atrial fiber architecture, first in idealized geometries \cite{pegolotti2019isogeometric} and then in realistic ones \cite{ferrer2015detailed,zygote2014}. 
We analysed the influence of atrial fiber bundles in the electric signal propagation by means of EP simulations performed on realistic geometries.  
Finally, we studied how a change in size of a single bundle affects the total activation sequence.

\begin{table}[t!]
	\begin{center}
		\begin{tabular}{ |c|c|c|c||c|c|c|c|c|c|c|c|c| } 
			\hline
			\bf{LA}  & $\textcolor{black}{\tau_{mv}}$ & $\tau_{lpv}$ & $\tau_{rpv}$ & \bf{RA} & $\textcolor{black}{\tau_{tv}}$ & $\tau_{icv} $ & $\tau_{scv}$ & $\tau_{ct_+}$ & $\tau_{ct_-}$ & $\textcolor{black}{\tau_{ib}}$ & $\textcolor{black}{\tau_{ras}}$ & $\textcolor{black}{\tau_{raw}}$  \\
			\hline
			\hline
			\bf{Ideal} & 0.65 & 0.65 & 0.10 & \bf{Ideal} & 0.90 & 0.90 & 0.10 & -0.10 & -0.18 & \textcolor{black}{0.35} & 0.135 & 0.55 \\  
			\hline
			\bf{Zygote}  & 0.85 & 0.85 & 0.20 & \bf{Zygote} & 0.90 & 0.85 & 0.30 & -0.55 & -0.60 & -0.25 & -0.10 & 0.60 \\
			\hline
			\textcolor{black}{\bf{Riunet}}  & \textcolor{black}{0.85} & \textcolor{black}{0.85} & \textcolor{black}{0.20} & \textcolor{black}{\bf{Riunet}} & \textcolor{black}{0.89} & \textcolor{black}{0.90} & \textcolor{black}{0.20} & \textcolor{black}{-0.10} & \textcolor{black}{-0.13} & \textcolor{black}{0.06} & \textcolor{black}{0.13} & \textcolor{black}{0.55} \\
			\hline
		\end{tabular}
	\end{center}
	\caption{Bundles parameters used for fibers generation in the idealized (Ideal) and realistic \textcolor{black}{(Zygote and Riunet)} LA and RA.}
	\label{table_atria_params}
\end{table}
\subsubsection{Atrial fibers generation}

We began applying the novel atrial LDRBM on idealized geometries. To construct them, we started by the surface representations of the right and left atrium generated as separated NURBS patches, as reported in \cite{pegolotti2019isogeometric,patelli2017isogeometric}. For each atrium, we created the corresponding triangular mesh using the constructive geometry module of \texttt{Netgen} (\url{https://ngsolve.org}). We considered this triangular mesh as the endocardium of our 3D model. To generate the \textcolor{black}{atrial} epicardium we extruded (using the \texttt{vmtk} software \cite{antiga2008vascular}) the endocardial surface by $2$ mm, which correspond to an average thickness of the atrial wall \cite{hoermann2019automatic,beinart2011left}. Finally, we  produced 3D tagged hexahedral mesh following the preprocessing pipeline described in Section \ref{sec:setting}, and then we applied our atrial LDRBM, see Figures \ref{fig:Atria_Fibers}(a-d). 

Afterwards, we treated the case of realistic left and right atria taken from the Zygote 3D heart model~\cite{zygote2014} \textcolor{black}{and from the Riunet repository (\url{https://riunet.upv.es/handle/10251/55150}). In particular, concerning the latter, we extracted both the endocardium and the epicardium (using \texttt{ParaView}), we removed all the inter-atrial connections (using \texttt{vmtk}) and then we created the corresponding 3D tagged hexahedral mesh (using the pipeline described in Section \ref{sec:setting}). Fibers generated by our atrial LDRBM are shown in Figures~\ref{fig:Atria_Fibers}(e-h) and \ref{fig:Atria_Fibers}(i-n) for the Zygote and Riunet geometries, respectively.}

The input values of the parameters $\tau_{i}$, which define the bundles dimension of the atrial LDRBM, are reported in Table \ref{table_atria_params}. We observe that the atrial LDRBM qualitatively capture the complex arrangement of fiber directions in almost all the principal anatomical atrial regions (see Figure~\ref{fig:Atria_Fibers}): for RA, superior (SCV) and inferior caval veins (ICV), tricuspid valve ring (TV), right appendage (RAA), septum (RAS), inter-caval bundle (IB), crista terminalis (CT), isthmus (IST) and atrial later wall (RAW); for LA, left (LPV) and right pulmonary veins (RPV), mitral valve ring (MV), left appendage (LAA), septum (LAS), septum wall (LSW), atrial lateral wall (LAW) and roof (LAR). 

\textcolor{black}{Following the atrial LDRBM rules defined at the end of Section \ref{sec:atria}}, circular fiber arrangements are exhibited around LPV, RPV, SCV, ICV, TV, MV, and encircle both appendages (RAA and LAA), see Figures \ref{fig:Atria_Fibers}(a-d) and \ref{fig:Atria_Fibers}(e-g). Fibers direction of CT runs longitudinally from the base of the SCV to the ICV, see Figures \ref{fig:Atria_Fibers}(d) and \ref{fig:Atria_Fibers}(h). RA structures like the IB and RAW are almost vertically oriented, whereas those of RAS are parallel to the CT, see Figures \ref{fig:Atria_Fibers}(c-d) and \ref{fig:Atria_Fibers}(g-h). IST fibers have the same direction of those of the TV, see Figures \ref{fig:Atria_Fibers}(d) and \ref{fig:Atria_Fibers}(h). The LAS fibers are aligned with the adjacent region of RAS, see Figures \ref{fig:Atria_Fibers}(a) and \ref{fig:Atria_Fibers}(e). Directions of LAR and LAW descend perpendicularly to MV (Figures \ref{fig:Atria_Fibers}(b) and \ref{fig:Atria_Fibers}(f)), while fibers of LSW present a smooth transition going to LAS and LAA (Figures~\ref{fig:Atria_Fibers}(a) and \ref{fig:Atria_Fibers}(e)). 

\subsubsection{Atrial electrophysiology} 
\begin{figure}[t]	
	\centering
	\includegraphics[width=0.95\textwidth]{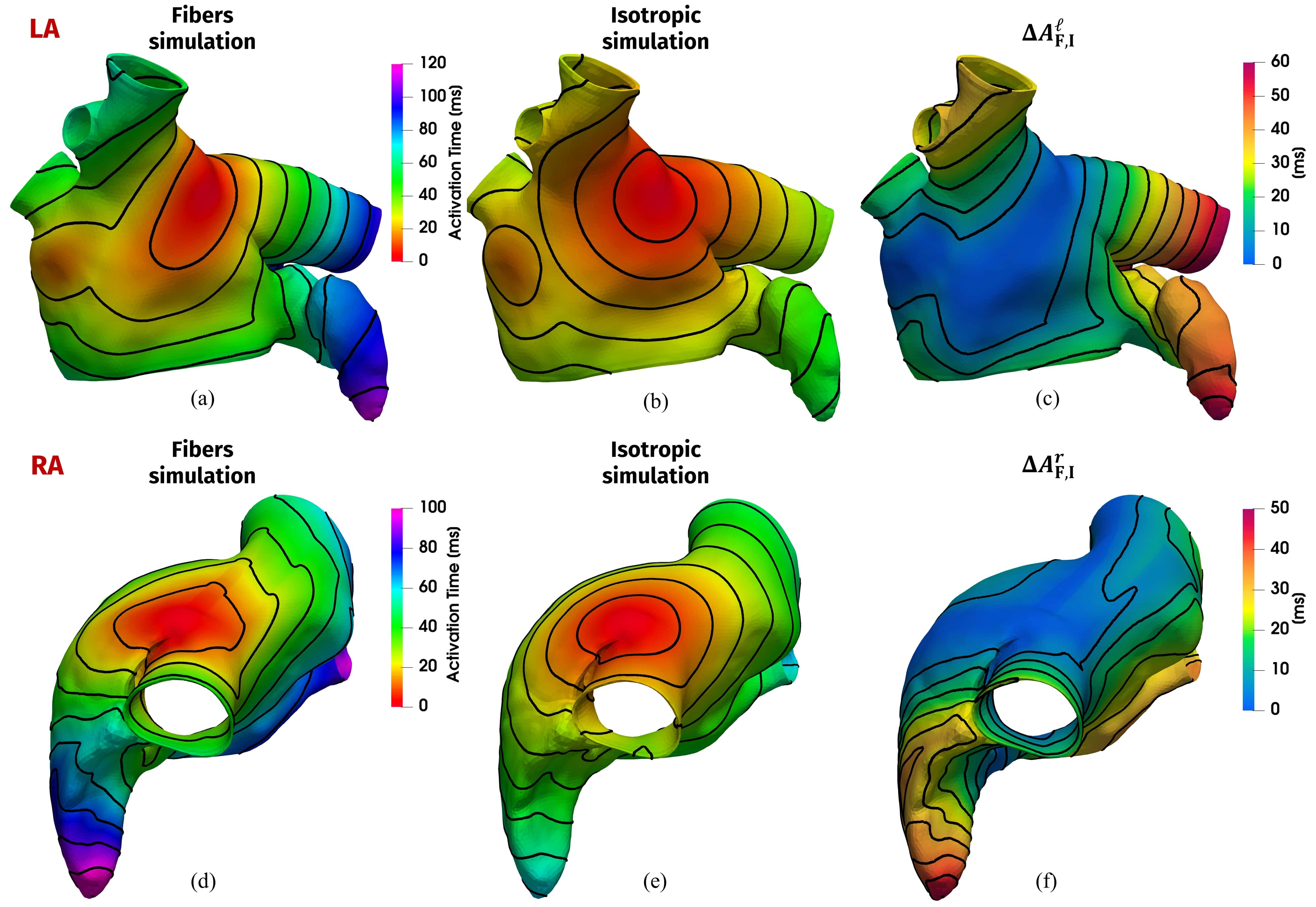}
	\caption{Comparison between the activation maps of EP simulations performed with the atrial LDRBM fiber generation and the isotropic model, \textcolor{black}{Zygote geometry}. Left (a,d): Fibers simulation; Centre (b,e): Isotropic simulation. Right (c,f): absolute difference $\Delta A^i_{\mathrm{F,I}}$ in the activation times for LA ($i=\ell$) and RA ($i=r$). Top (a,b,c): LA; Bottom (d,e,f): RA.}
	\label{fig:Atria_AT}
\end{figure}
\begin{figure}[t!]	
	\centering
	\includegraphics[width=0.95\textwidth]{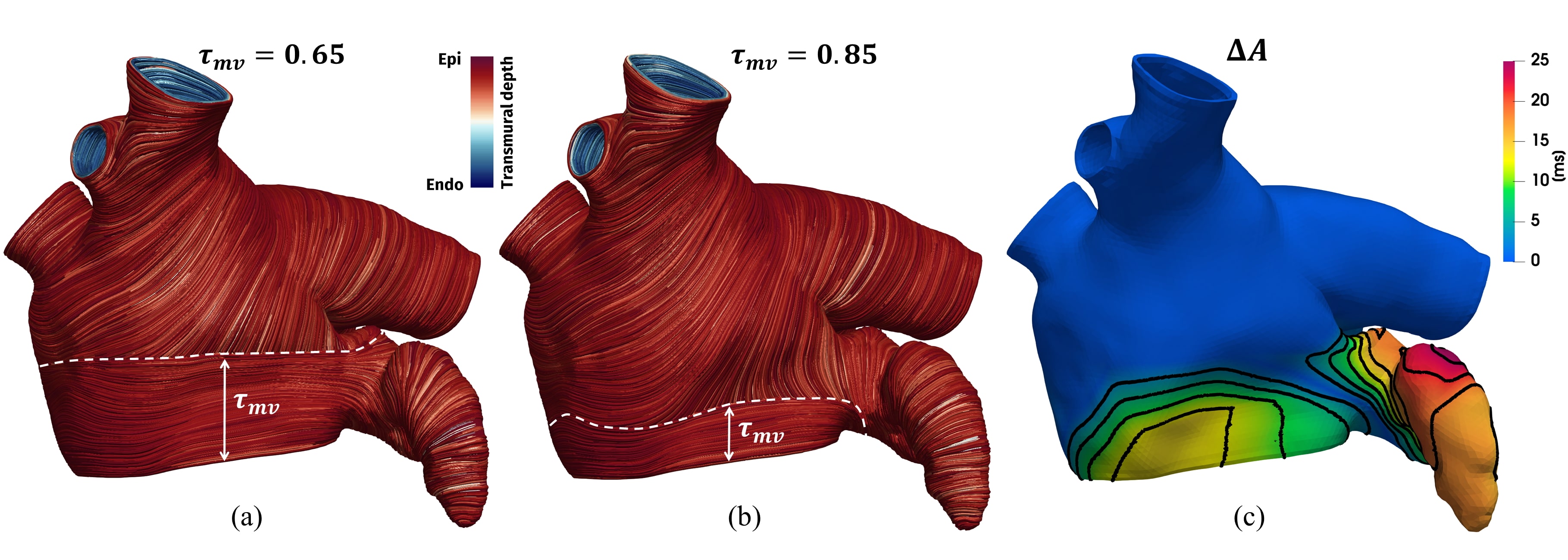}
	\caption{Comparison between EP simulations with different values of \textcolor{black}{$\tau_{mv}$} in the atrial LDRBM fiber generation for the Zygote LA. \textcolor{black}{Dashed lines in white represent the limit of the MV bundle.} Left (a): \textcolor{black}{$\tau_{mv}$}$=0.65$; Centre (b): \textcolor{black}{$\tau_{mv}$}$=0.85$; Right (c): absolute difference $\Delta A$ in the activation times.}
	\label{fig:Atria_Tau}
\end{figure}  
In order to analyse the influence of atrial fiber bundles in the electric signal propagation we performed several EP simulations (with the setting specified in Section \ref{sec:setting}) on the realistic Zygote atrial geometries. 

Firstly, we made a comparison with an isotropic model. For the atrial LDRBM, we considered the parameters detailed in Table \ref{table_atria_params}, while the isotropic simulations were carried out by setting in \eqref{D_tensor} $\sigma_f=\sigma_s=\sigma_n=7.0$ mS/cm, that is a representative value chosen for the conductivity along the atrial fiber direction (see Table \ref{table_conductivity} in the Appendix). To initiate the signal propagation in RA we applied a single stimulus in the Sino-Atrial-Node (SAN) which lies in the musculature of CT at the anterolateral junction with the SCV~\cite{ho2009importance}. For LA we stimulated the main inter-atrial connections: the Bachman's Bundle (BB), located in the LSW; the upper part of the Fossa Ovalis (FO) in the centre of LAS; the Coronary Sinus Musculature limbs (CSM) placed at the bottom of LAW \cite{sakamoto2005interatrial}. Activation of FO and CSM was delayed, with respect to the BB stimulus, by 14 ms and 52 ms, respectively. 

Figure \ref{fig:Atria_AT} displays the results of the comparison among simulations performed with the atrial LDRBM fibers and the isotropic model for both RA and LA. Both the activation pattern and activation time present significant differences.
To provide a quantification, we computed the absolute difference $\Delta A^i_{\mathrm{F,I}}$ in the activation time:
\begin{equation}
\Delta A^i_{\mathrm{F,I}}(\boldsymbol{x})=|A^i_{\mathrm{F}}(\boldsymbol{x})-A^i_{\mathrm{I}}(\boldsymbol{x})| \qquad i=r,\ell, 
\end{equation} 
where $i=r,\ell$ refer to LA ($i=\ell$) and RA ($i=r$) and $A^i_{\mathrm{F}}$ and $A^i_{\mathrm{I}}$ are the activation times obtained by the simulations with and without fibers, respectively. Most of the differences occur at LPV and LAA for LA, and at RAA and TV for RA. Finally, we computed the maximal discrepancy, $M^i_{F,I}=~\max_{\boldsymbol{x}\in \Omega_{atrial}}\Delta A^i_{F,I}(\boldsymbol{x}),\,i=~\ell,r$:
$$M^\ell_{F,I}=60\, \mathrm{ms} \,(52 \, \%), \qquad M^r_{F,I}=48 \, \mathrm{ms} \,(44 \, \%),$$
where in brackets we reported the relative values computed as $M^i_{F,I}/A^i_{max}$, with $A^\ell_{max}=116$ ms and $A^r_{max}=108$ ms the total activation times. For RA  $M^r_{F,I}$ is placed in RAA, while for LA $M^\ell_{F,I}$ is located in LPV.  

Then, we investigated how a local change in a single LA bundle (the \textcolor{black}{mitral valve} one) affects the total activation pattern. We performed two EP simulations with the same fiber setting used for the comparison with an isotropic model, except for the value of \textcolor{black}{$\tau_{mv}$}, which was set equal to $0.65$ and $0.85$. Figures~\ref{fig:Atria_Tau}(a-~b) depict the corresponding generated fibers: notice that with \textcolor{black}{$\tau_{mv}$}$=0.65$ the \textcolor{black}{MV} bundle is thicker with respect to the one obtained with \textcolor{black}{$\tau_{mv}$}$=0.85$, see Figures \ref{fig:Atria_Tau}(a-b). We also reported the absolute difference in the computed activation times for the two fiber architectures, see Figure \ref{fig:Atria_Tau}(c). The maximal discrepancy, located in LAA, is of $28$~ms which corresponds to $24\%$ of the total activation time for LA ($116$ ms).

\begin{figure}[t!]	
	\centering
	\includegraphics[width=1\textwidth]{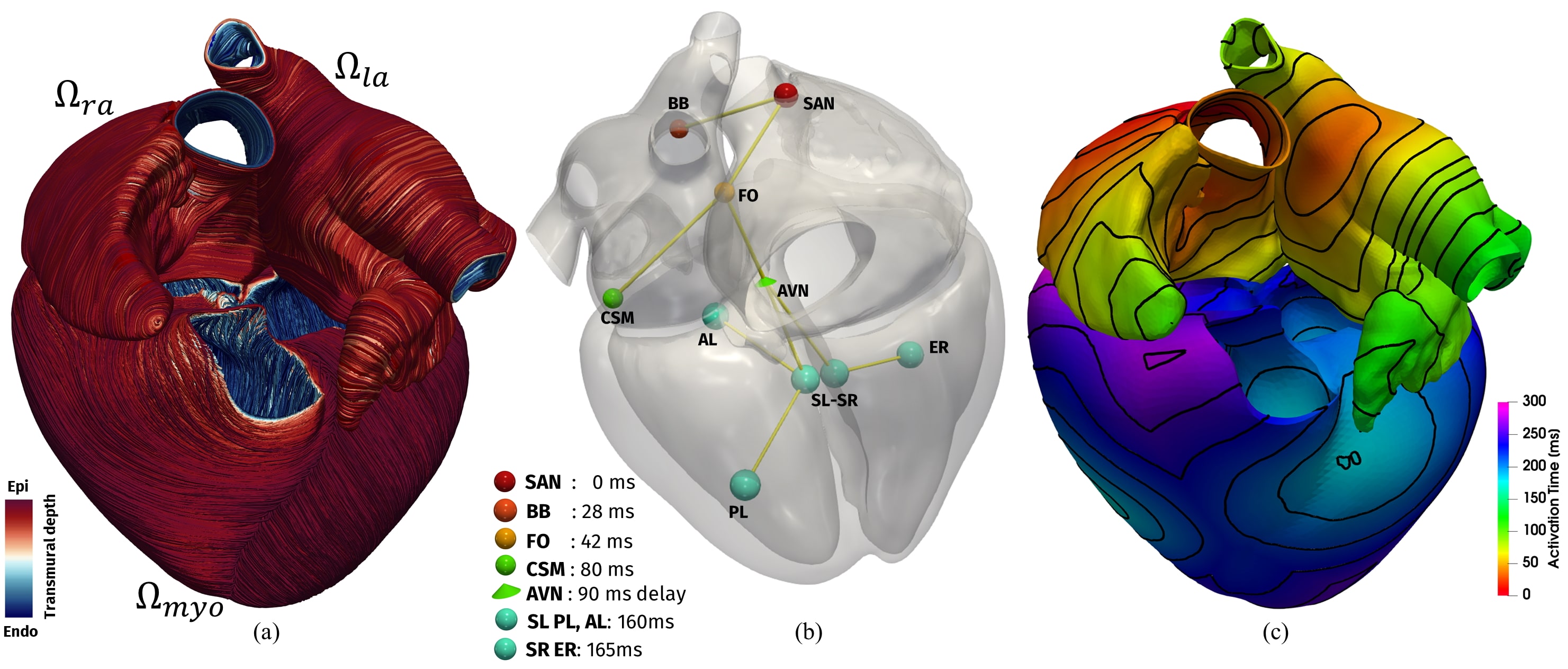}
	\caption{Left (a): Fiber generation applied to realistic Zygote Heart model; LDRBM by Doste (D-RBM) was applied to reconstruct the ventricular fibers on $\Omega_{myo}$, while the atrial LDRBM was employed for the fiber architecture on $\Omega_{ra}$ and $\Omega_{la}$. Centre (b): Stimuli applied in the four chamber model to mimic the Cardiac Conduction System (CCS) pathway; SAN: Sino-Atrial Node; BB: Bachmann’s Bundle; FO: Fossa Ovalis; CSM: Coronary Sinus Musculature; AVN: Atrio-Ventricular Node; AL: Left Anterior; PL: Left Posterior; SL, SR: Left and Right Septum; ER: Right Endocardium. Right (c): activation maps computed from EP simulation.}
	\label{fig:Heart_Fibers}
\end{figure} 
\subsection{Whole heart fibers and electrophysiology}\label{sec:heart}
In this section we present the whole heart fiber generation, using LDRBMs for both atria and ventricles. Moreover, we show an EP simulation using physiological activation sites and including the fiber generated by LDRBMs. We use the Zygote heart model \cite{zygote2014} both for the full biventricular and the \textcolor{black}{atrial} geometries introduced in Sections \ref{sec:real-biv} and \ref{sec:atria}, respectively. 

\subsubsection{Whole heart fibers}
For the fibers generation, we adopted D-RBM \cite{doste2019rule} (see Section \ref{sec:ventr}) for the ventricles, with the same setting of Section \ref{sec:real-biv}, and the proposed
LDRBM for the atria (see Section \ref{sec:atrial}), with the configuration of Section~\ref{sec:atria} (see also Table \ref{table_atria_params}). Figure \ref{fig:Heart_Fibers}(a) displays the heart geometry equipped with the prescribed LDRBMs fibers. 
\begin{figure}[t!]	
	\centering
	\includegraphics[width=1\textwidth]{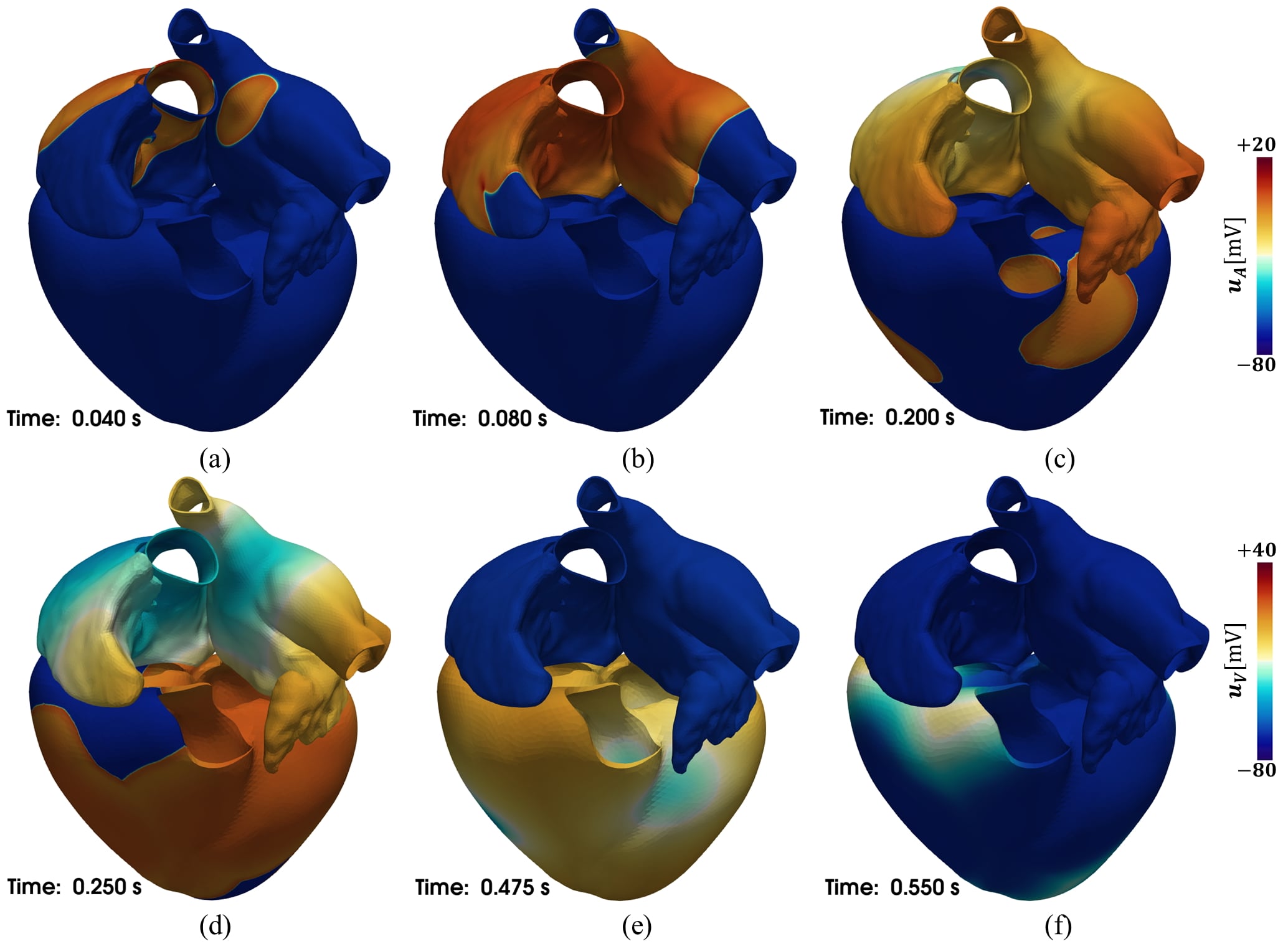}
	\caption{Evolution of the transmembrane potentials for the ventricles $u_V$ and for the atria $u_A$ in the Zygote heart model during a cardiac cycle.}
	\label{fig:Heart_u}
\end{figure} 
\subsubsection{Whole heart electrophysiology}
To model the EP activity in the cardiac tissue we used the monodomain equation endowed with the TTP and CRN ionic models for the ventricles and for the atria, respectively, with the settings described in Section~\ref{sec:setting}.          

The interactions among atria and ventricles are based on the following assumptions on the cardiac conductions system (CCS) connections, showed in Figure \ref{fig:Heart_Fibers}(b). The ventricles are electrically isolated from the atria by the atrioventricular grooves \cite{anderson2000anatomy}; the atria are electrically separated by the insulating nature of the atrial septum (dividing RA from LA) apart from muscular continuity at the rim of Fossa Ovalis~\cite{ho2009importance}. The CCS pathway was modelled as a series of spherical delayed stimuli along the heart geometry that mimic the inter-atrial connections, the \textcolor{black}{Atrio-Ventricular Node (AVN)} delay and the main area of ventricular electrical activation: specifically, when the transmembrane potential front reaches these points a stimulus current is triggered, see Figure \ref{fig:Heart_Fibers}(b). \textcolor{black}{Notice that, although the Purkinje network system should be used to account for a detailed ventricular activation (see e.g. \cite{VerPalCat14,PasRomSeb1110}), for the sake of simplicity several endocardial stimuli were applied to each ventricle \cite{durrer1970total}.}

The CCS electric signal originates at the Sino-Atrial Node (SAN, $t=0$ ms) and travels from RA to LA through three inter-atrial connections, the Bachmann’s Bundle (BB, $t=28$ ms), the rim of Fossa Ovalis (FO, $t=42$ ms) and the Coronary Sinus Musculature (CSM, $t=80$ ms) \cite{ho2009importance,sakamoto2005interatrial}. When the electric signal reaches the AVN, located at the lower back section of the inter-atrial septum near the coronary sinus opening, it is subject to a delay ($90$~ms) that allows the complete activation of the atria before ventricles electric propagation starts~\cite{ho2009importance}. Finally, \textcolor{black}{ventricular} endocardial areas are activated: in the anterior para-septal wall (AL), in the left surface of inter-ventricular septum (SL) and in the bottom of postero-basal area (PL), for the left ventricle ($t=160$ ms); in the septum (SR) and in the free endocardial wall (ER), for the right ventricle ($t=165$ ms) \cite{durrer1970total}.

Figure \ref{fig:Heart_Fibers}(c) depicts the activation maps computed by the heart EP simulation. The complete \textcolor{black}{atrial} depolarization occurs after about $120$ ms, while that of ventricles after about $t=270$ ms. The last region to be activated is LAA for the atria, while the postero-basal area of the right ventricle for the ventricles. 

The transmembrane potentials evolution for the ventricles $u_V$ and for the atria $u_A$ are shown in Figure~\ref{fig:Heart_u}. As expected, the electric signal initiates at the SAN and spreads from right to left atrium, see Figures \ref{fig:Heart_u}(a) and \ref{fig:Heart_u}(b). Then, after the delay at the AVN, the ventricles start to activate, see Figure \ref{fig:Heart_u}(c). The atrial repolarization arises during ventricular depolarization, see Figures \ref{fig:Heart_u}(c-e). Finally, after the isoelectric ventricular activity, the whole heart return to the depolarized initial configuration, see Figure \ref{fig:Heart_u}(f).     

\section{\textcolor{black}{Discussion}}\label{sec:disc} 

\subsection{\textcolor{black}{Comparison among the three ventricular LDRBM}}
Three existing LDRBMs for fibers generation in the ventricles (R-RBM, B-RBM, D-RBM) were reviewed by means of a communal mathematical description. Some extensions were introduced allowing the inclusion of different fiber orientations in the left and right ventricles for R-RBM and B-RBM, the rotation of all the myofiber vectors for R-RBM, and the fibers generation up to cardiac valve rings for B-RBM.
 
The comparison among the three ventricular LDRBMs showed that R-RBM and B-RBM were able to recover almost the same fiber orientations of D-RBM thanks to our extensions (Figures \ref{fig:BiV_Ideal_Fibers} and \ref{fig:BiV_Zygote_Fibers}). However, some local differences persist in the methods utilized. Specifically, most of the fibers orientation discrepancies were evident in the right ventricular endocardium facing to the septum, in the inter-ventricular junctions, in the right epicardial lower region and in the proximity of the valve rings (Figures 7 and 9).

To better discuss the differences revealed at the ventricular septum, the main specific issues characterizing the methods were: i) the presence (or absence) of a different septal fiber orientations in the left and right ventricles and ii) the presence (or absence) of a fibers discontinuity in the inter-ventricular septum. These features were determined by the different definitions of transmural and normal directions implemented by the three LDRBMs (Section 2). Indeed, by construction, R-RBM fibers in the whole septum belong to the left ventricle, while for B-RBM and D-RBM they are shared between the two ventricles. Consequently, according with the anatomical observations \cite{scollan2000reconstruction,ho2006anatomy,sanchez2015anatomical,lunkenheimer2013models,stephenson2016functional}, we argue that the latter two LDRBMs better describe the fibers direction in the right endocardium (Figures \ref{fig:BiV_Ideal_Fibers} and \ref{fig:BiV_Zygote_Fibers}).  Moreover, unlike the other LDRBMs producing a smooth transition in the fiber field around the inter-ventricular junctions, D-RBM showed a strong discontinuity in the transition across the two ventricles (Figures~\ref{fig:BiV_Ideal_Fibers} and \ref{fig:BiV_Zygote_Fibers}). Standard anatomical observations claim that the fibers are almost continuous through the septum \cite{ho2006anatomy,sanchez2015anatomical}, though recent studies support the thesis of septal fibers discontinuity \cite{kocica2006helical,boettler2005new}. Further investigations should be performed to establish which LDRBM better reproduces the anatomical ground truth.

The different orientations of the septal fibers produced a significant impact in terms of EP numerical results. A dissimilar activation pattern and timing were observed at the septum, especially in the right endocardial region, confirmed by the highest discrepancy between R-RBM and the other two methods reaching 28-29\% of the total activation time (Figure \ref{fig:BiV_Ideal_AT}). Conversely, D-RBM and B-RBM yield almost the same activation pattern with only 10\% timing discrepancy, thanks to the extensions introduced in B-RBM (Figures \ref{fig:BiV_Ideal_AT} and~\ref{fig:BiV_Zygote_Fibers}). All the above results confirmed the importance of including a specific fiber orientation in the right ventricle with respect to the fiber orientation in the left ventricle. 

A different septal activation could have a strong influence also in the mechanical contraction triggered by the electrical activation of the myocardium. Therefore, the importance of including suitable fiber orientations for the right ventricle, together with the discontinuity around the inter-ventricular junctions and specific orientations in the ventricular outflow tract should be investigated also in the context of the electromechanical simulations.

\subsection{\textcolor{black}{Comparison of atrial fibers with pre-existing results and anatomical data}}
\begin{figure}[t!]	
	\centering
	\includegraphics[width=0.95\textwidth]{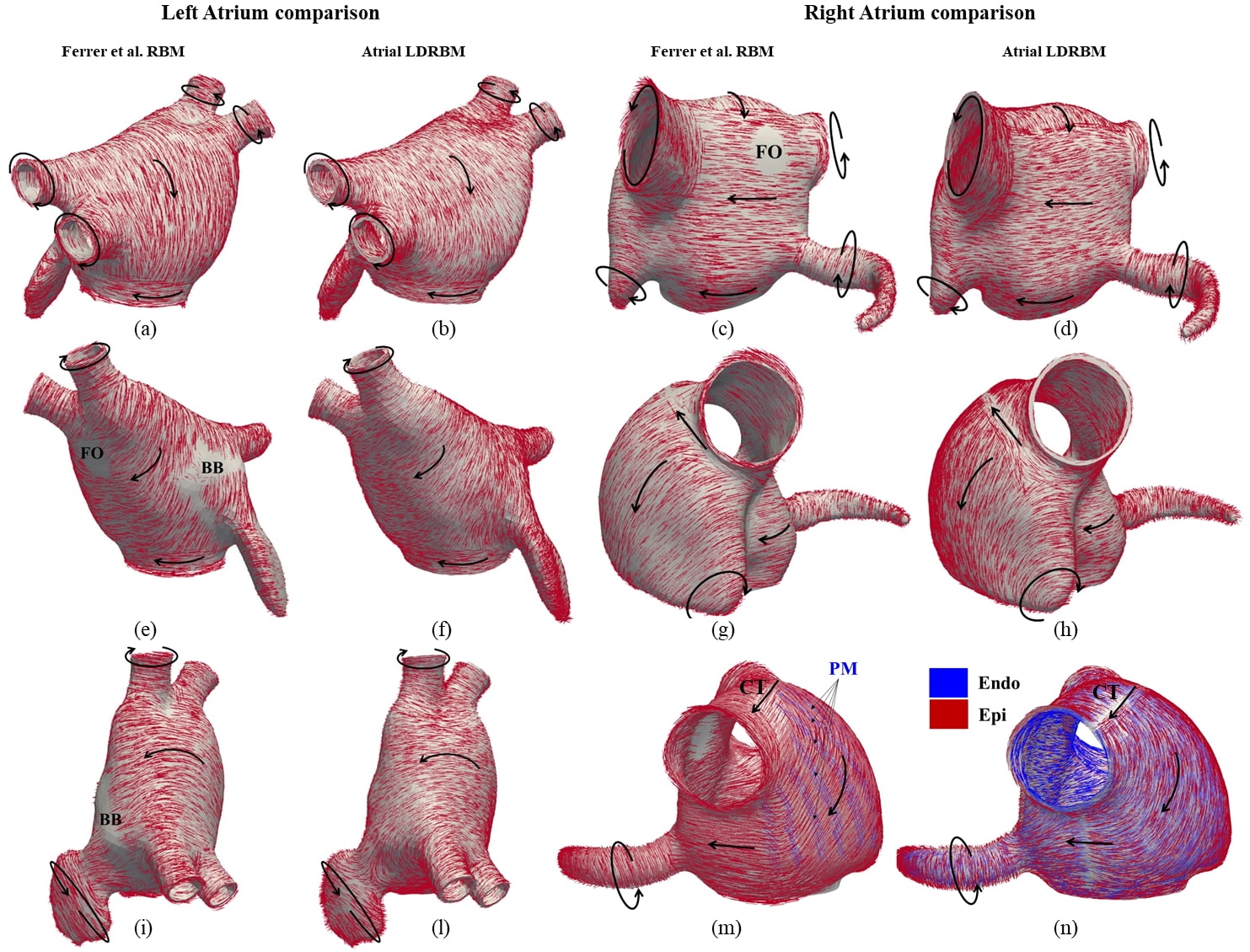}
	\caption{\textcolor{black}{Fiber orientations (rendered in red glyphs) comparison in the Riunet geometry between the atrial LDRBM and RBM in \cite{ferrer2015detailed}. First column (a,e,i): RBM in \cite{ferrer2015detailed} left atrium; Second column (b,f,l): atrial LDRBM left atrium; Third column (c,g,m): RBM in \cite{ferrer2015detailed} right atrium; Fourth column (d,h,n): atrial LDRBM right atrium; Black arrows represent the main fibers direction in specific atrial regions; FO: Fossa Ovalis; BB: Bachmann's Bundle; CT: Crista Terminalis; PM: Pectinate Muscles.}}
	\label{fig:Confronto}
\end{figure}
We presented a novel atrial LDRBM where suitable intra-atrial distance functions are introduced by means of harmonic problems. To represent the fiber bundles, the novel methodology used the gradient of the intra-atrial distances combined with a precise definition of the boundary sections where boundary conditions are prescribed for the harmonic problems. This strategy makes the fibers less open to subjective variability. On the other hand, the bundles dimension could be adapted case by case changing the parameters involved in Algorithms 1 and 2 (see Section \ref{sec:atrial}). Therefore, unlike most of the existing atrial RBM requiring manual or semi-automatic interventions, our new method can be easily applied to any arbitrary geometry.

In order to verify the reliability of the atrial LDRBM, the results obtained by our novel fiber generation strategy in the Riunet geometry were compared with those previously reported in \cite{ferrer2015detailed} \footnote{Freely available online at \url{https://riunet.upv.es/handle/10251/55150}} with a RBM incorporating a detailed regional description of fiber directions provided by the anatomical observations~\cite{tobon2013three,ferrer2015detailed}. Specifically, the RBM in \cite{ferrer2015detailed} was built by a manual subdivision of the atrial geometry in several regions to embed a detailed fibers description \cite{ferrer2015detailed}. This model, unlike the presented atrial LDRBM, includes also the inter-atrial conduction bundles, i.e. the Bachmann's Bundle (BB), the limb of Fossa Ovalis (FO) and the connections of the Coronary Sinus Musculature (CSM). To account for the presence of such conduction bundles, in the present work their effect were surrogated by using a suitable timed pacing of LA \cite{pegolotti2019isogeometric}. This artificial atrial conduction system is suitable to generate a physiological activation at least in case of sinus rhythm \cite{pegolotti2019isogeometric}.

The fibers field generated by our LDRBM (Figure \ref{fig:Confronto}), owing to a suitable choice of the input parameters $\tau_{i}$ (Table \ref{table_atria_params}), is in excellent agreement with the finding of the RBM previously proposed \cite{ferrer2015detailed}, reproducing almost the same fiber orientations among the different atrial bundles. Major differences are visible in a restricted zone of BB, embracing the left atrial appendage (Figures~\ref{fig:Confronto}(e,f,i,l)), at the limb of FO (Figures~\ref{fig:Confronto}(c,d,e,f)), and in the right atrial wall endocardium due to the Pectinate Muscles (PM) which, differently from the present work, were explicitly included in \cite{ferrer2015detailed} (Figures~\ref{fig:Confronto}(m,n)). As for the PM, which should be roughly perpendicular to CT, we surrogated them by prescribing, through Algorithm 1, the fibers directions in the right atrial wall endocardium to be perpendicular to CT.

\begin{figure}[t!]
	\centering
	\includegraphics[width=0.85\textwidth]{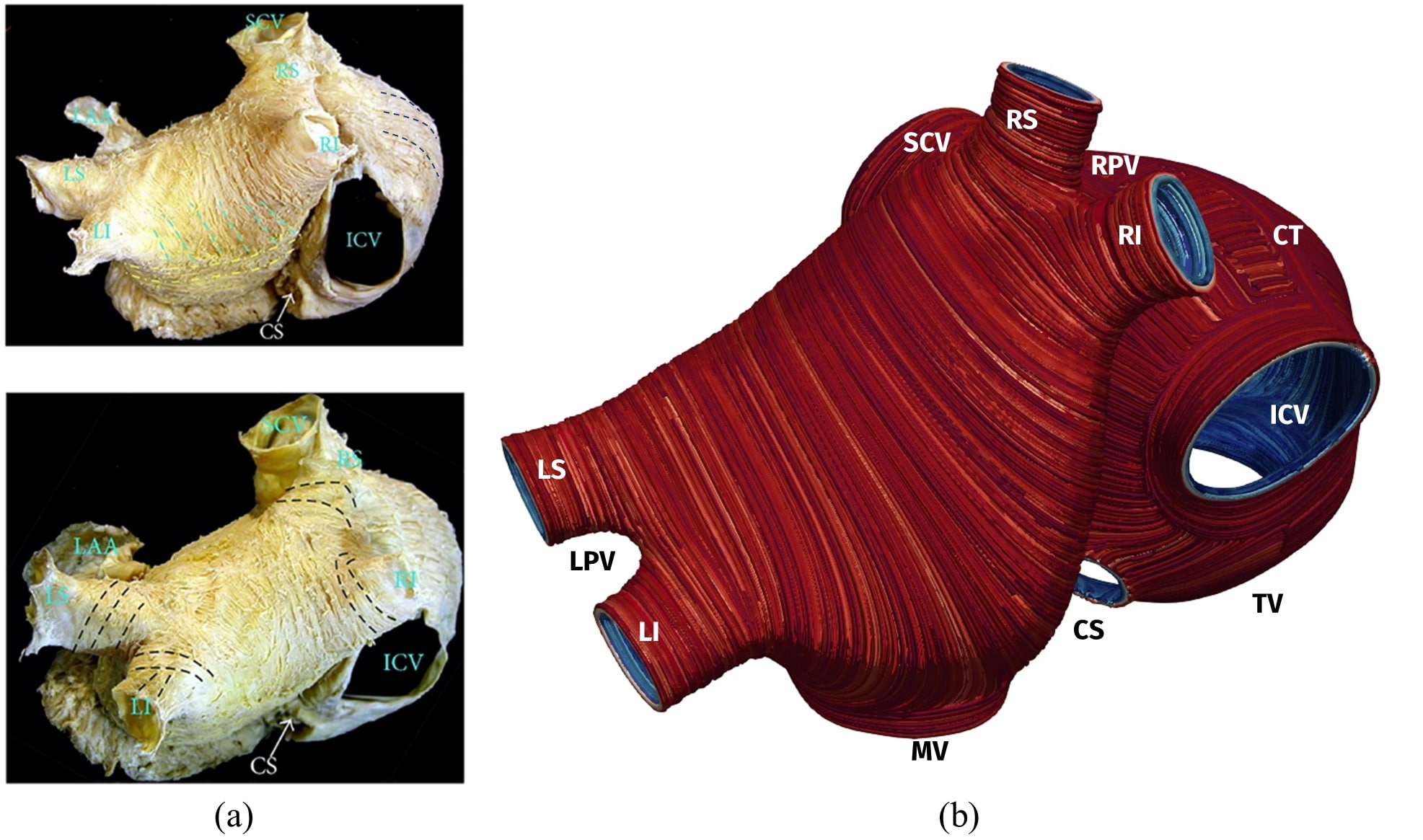}
	\caption{\textcolor{black}{Comparison between anatomical atrial dissections in a normal human atria and fibers obtained by the atrial LDRBM in the Riunet geometry: (a) atrial fiber structure in the normal human heart where dashed lines display the main fibers orientations; (b) fiber prescribed by the atrial LDRBM rules (R1-R6). ICV, SCV: inferior and superior caval veins; CT: crista terminalis; CS: coronary sinus; TV, MV: tricuspid and mitral valve rings; LPV: left pulmonary superior (LS) and inferior (LI) veins; RPV: right pulmonary superior (RS) and inferior (RI) veins. Figures on the left are from the work by Sánchez‐Quintana \cite{sanchez2014left}.}}
	\label{fig:Histo_F}
\end{figure}

In addition, a graphical comparison of our results obtained in the Riunet geometry was provided with
anatomical pictures of atrial dissections in a normal human heart taken from \cite{sanchez2014left} (Figure \ref{fig:Histo_F}). The fiber directions predicted by the atrial LDRBM rules showed an excellent agreement with the anatomical study, sustaining the validity of rules R1-R6 (Section \ref{sec:atria}).

The atrial EP simulations revealed a strong influence of the complex atrial fiber architecture on the
electric signal propagation. The activation pattern and timing, induced by the atrial LDRBM fibers, were
consistently different from the isotropic model (Figure \ref{fig:Atria_AT}) in accordance with the previous findings~\cite{krueger2011modeling,fastl2018personalized}. Electrophysiology simulations, embedded with the atrial LDRBM fibers, provided a total activation time of 108 ms for RA and 116 ms for LA (Figure \ref{fig:Atria_AT}). These values are compatible with the timings
predicted in~\cite{collin2013surface} for RA (102 ms) and in \cite{krueger2011modeling} for LA (115 ms).

\subsection{\textcolor{black}{The whole heart case}}

Numerical results including the full heart LDRBMs fiber generations and a related EP simulation with physiological activation sites were illustrated, producing a physiologically compatible timing for the heart activation \cite{quarteroni2017integrated,franzone2014mathematical,quarteroni2019}. The complete atrial depolarization occured after about 120 ms, while the ventricular activation after about 270 ms (Figure \ref{fig:Heart_Fibers}). The last region activated was LAA for the atria, while the postero-basal area of the right ventricle for the ventricles (Figure \ref{fig:Heart_u}), both in accordance with previous reports \cite{sakamoto2005interatrial,durrer1970total}.

The proposed whole heart fibers generation methodology is computationally inexpensive, efficient and easy to implement, and it allows to include realistic cardiac muscle fibers architecture on whole heart geometries of arbitrary shape. As a consequence, it is possible to generate patient cohorts heart fibers, fed by input parameters inferred from histology or DTI studies, through an automated and computationally efficient pipeline. Our proposed methodology provides an important contribution to incorporate patient specific whole heart fiber field into EP and electromechanics simulations, allowing both the study of clinical cases as well as investigating medical questions.

\subsection{\textcolor{black}{Limitations}}\label{sec:limit}

This work, mainly focused on providing a modeling pipeline for the fibers generation in the whole heart,
presents the following limitations.

The presented version of the atrial LDRBM was developed to represent all prominent fiber bundles,
but some simplifications were introduced. In particular, LA and RA were treated as two independent entities and the inter-atrial connections (BB, FO, CSM) were not explicitly represented \cite{ferrer2015detailed,kruger2013personalized}, but their effect was rather
surrogated. Another limitation consisted in the use of a surrogate model for the PM instead of an
explicit modelization of such bundle. Finally, we were aware that in the presented atrial LDRBM framework it was possible, in principle, to prescribe a different transmural orientation in some bundles because we defined a transmural gradient, but this step was not considered to simplify the entire procedure. Nevertheless, this omission could be justified by the observation that the inclusion of a variable transmural orientation should not provide significant changes in EP simulations, as already highlighted \cite{fastl2018personalized}. All the above issues are currently under development to further improve the quality and validity of the present study.

The EP at cellular and tissue levels was modelled as homogeneous through the whole heart in order to
better focus and study the effect of the heart fiber architecture. Although in principle this assumption is a well accepted approximation for the ventricular tissues of healthy individuals, several atrial regions exhibit distinct EP properties. These should be taken into consideration to achieve a more realistic activation/repolarization pattern \cite{ferrer2015detailed,lemery2007normal}, particularly when reproducing pathological conditions affecting the atrial chambers \cite{dossel2012computational}.

\section{Conclusions}\label{sec:concl}
In this work, we provided a unified framework for generating cardiac muscle fibers in a full heart computational domain. This allowed us to obtain physically meaningful EP simulations of a four chambers heart realistic domain.

\textcolor{black}{We reviewed existing ventricular LDRBMs providing a communal mathematical description and introducing also some modeling improvements 
regarding the right ventricle.} We then carry out systematic comparisons of ventricular LDRBMs by means of numerical EP simulations, highlighting the importance of including a proper fiber orientation for the right ventricle. 

We proposed a novel LDRBM to be used for generating atrial fibers and we applied it to both idealized and realistic geometries showing, \textcolor{black}{through comparisons with another RBM and with anatomical atrial dissections}, that the atrial LDRBM capture the complex arrangement of fiber directions in almost all the anatomical atrial regions. We analysed the influence of \textcolor{black}{atrial} fiber bundles by means of EP simulations in a realistic geometry, verifying the strong effect of their complex architecture in the electric signal propagation.

Finally, we presented an EP simulation of a realistic four chamber heart including fibers generated by LDRBMs for both atria and ventricles. 
\begin{figure}[t!]	
	\centering
	\includegraphics[width=0.75\textwidth]{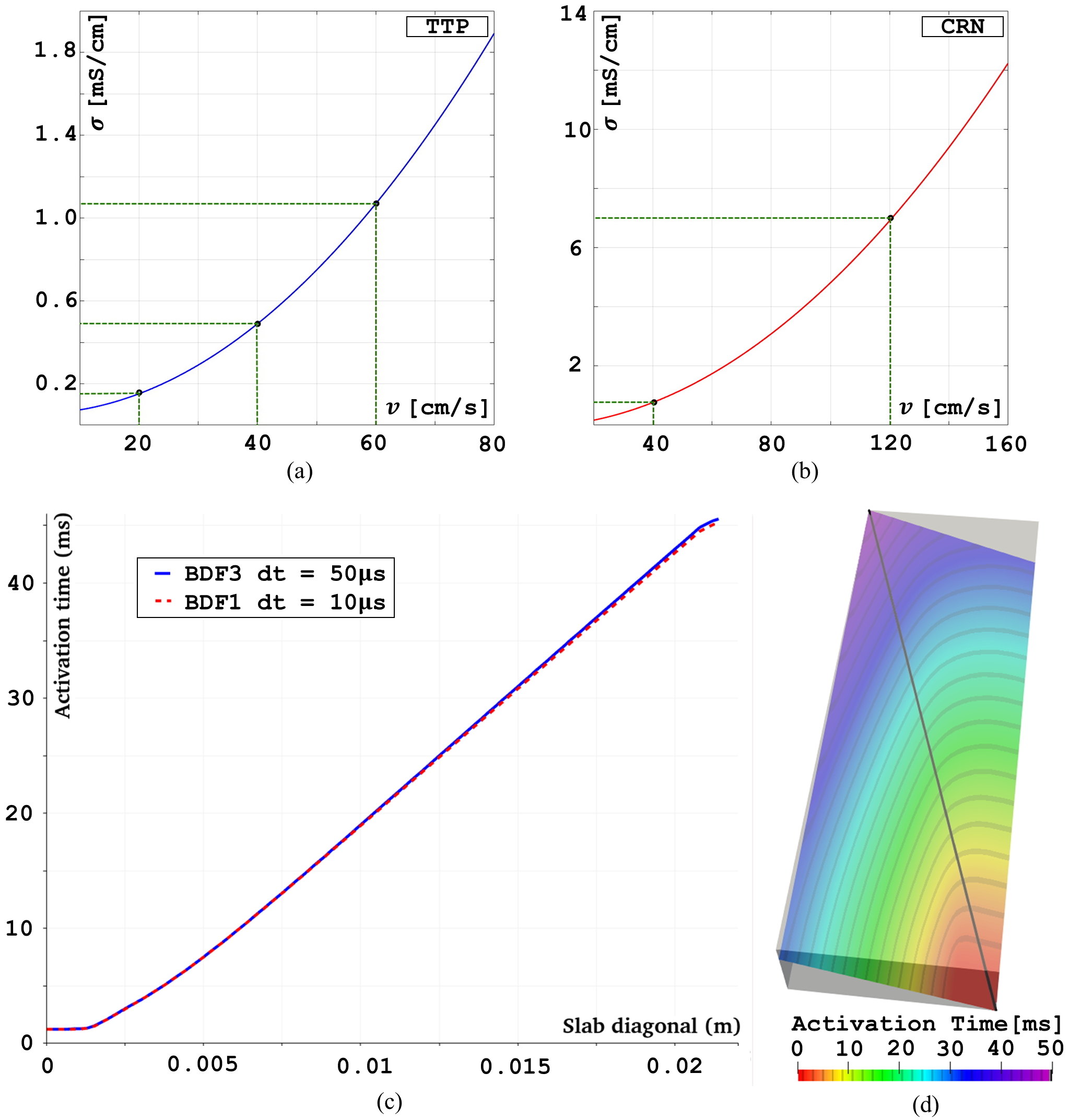}
	\caption{Top (a,b): Fitting procedure used to estimate the conductivity $\sigma$ required to match specific conduction velocity~$v$~\cite{costa2013automatic}; (a): using the TTP ionic model to obtain $60$, $40$ and $20$ cm/s; (b): using the CRN ionic model to obtain $120$ and $40$ cm/s. The values for $\sigma_f$, $\sigma_s$ and $\sigma_n$ are reported in Table \ref{table_conductivity}. Bottom (c,d): Comparison between $\text{BDF}3$ and $\text{BDF}1$ time discretization for the monodomain system \eqref{eq:monodomain_system} in the slab benchmark problem \cite{niederer2011verification}; (c): plot of the activation time alongside the slab diagonal (displayed in black on the Right); Red: $\text{BDF}1$; Blue: $\text{BDF}3$. (d): activation time in a clipped slice of the slab for $\text{BDF}3$ time discretization.}
	\label{fig:CV_estimate}
\end{figure}
\begin{table}[ht!]
	\begin{center}
		\begin{tabular}{ |c|c|c|c| } 
			\hline
			Muscle type (ionic model) & $\sigma_f$[mS/cm] & $\sigma_s$[mS/cm] & $\sigma_n$[mS/cm] \\
			\hline
			\hline
			Ventricles (TTP) & 1.07 & 0.49 & 0.16 \\
			\hline              
			Atria (CRN)        & 7.00 & 0.77 & 0.77 \\ 
			\hline					
		\end{tabular}
	\end{center}
	\caption{Conductivity values $\sigma_f$, $\sigma_s$ and $\sigma_n$ obtained after the fitting procedure, displayed in Figure \ref{fig:CV_estimate}, for the ventricles (using TTP) and for the atria (using CRN).}
	\label{table_conductivity}
\end{table}


\section*{\textcolor{black}{Appendix}}

\subsection*{\textcolor{black}{On the choice of physical parameters and numerical settings}}\label{sec:param}

The conductivity values $\sigma_f$, $\sigma_s$ and $\sigma_n$ were fitted by an iterative procedure described in \cite{costa2013automatic} (see also \cite{fastl2018personalized,augustin2016anatomically}) in order to match the following conduction velocity values: for the ventricles, $60$ cm/s in the fiber direction $\boldsymbol f$,  $40$ cm/s in the sheet direction $\boldsymbol s$ and $20$ cm/s in the normal direction $\boldsymbol n$ \cite{augustin2016anatomically}; for the atria, $120$~cm/s in the fiber direction $\boldsymbol f$ and $40$ cm/s along the sheet $\boldsymbol s$ and cross-fiber directions $\boldsymbol n$ \cite{fastl2018personalized,augustin2019impact,dimitri2012atrial}. In Figures \ref{fig:CV_estimate}(a) and \ref{fig:CV_estimate}(b) we show the results of this fitting procedure. The estimated values for $\sigma_f$, $\sigma_s$ and $\sigma_n$ are reported in Table \ref{table_conductivity}.

Regarding the mesh element size $h$ and the time step $\Delta t$, related to the
space and time discretizations of the system \eqref{eq:monodomain_system}, accuracy constraints are imposed when biophysical models (as CRN \cite{courtemanche1998ionic} and TTP \cite{ten2006alternans}) are used: $h=100$--$500$ $\mu$m and $\Delta t=1$--$50$ $\mu$s \cite{plank2008mitochondrial,vigmond2008solvers,niederer2011verification}. These strong restrictions are motivated mainly by the fast upstroke of cellular depolarization which produces a step-like wavefront over a small spatial extent~\cite{lopez2015three}. For the space discretization, we used continuous FE of order 1 ($Q_1$) on hexahedral meshes with an average mesh size of $h=350$ $\mu$m, an acceptable value at least for linear finite element approximation and for physiological cases \cite{trayanova2011whole,quarteroni2017integrated,hurtado2018non,jilberto2018semi,arevalo2016arrhythmia}. Concerning the time discretization, we used the $\text{BDF}$ of order $\sigma=3$ with a time step of $\Delta t=50$ $\mu$s. Although, the most common time discretization used in literature for the monodomain system \eqref{eq:monodomain_system} is the $\text{BDF}1$ (commonly known as the backward Euler method), which requires a time step at most of $10$ $\mu$s \cite{vigmond2008solvers}, BDF3 allows us to use a larger value of $\Delta t$ to obtain the same accuracy of $\text{BDF}1$. To confirm this, in Figures \ref{fig:CV_estimate}(c) and \ref{fig:CV_estimate}(d) we report a comparison between BDF3 with $\Delta t=50$~$\mu$s and BDF1 $\Delta t=10$ $\mu$s on a benchmark problem proposed in \cite{niederer2011verification}.    

\begin{figure}[t!]	
	\centering
	\includegraphics[width=0.75\textwidth]{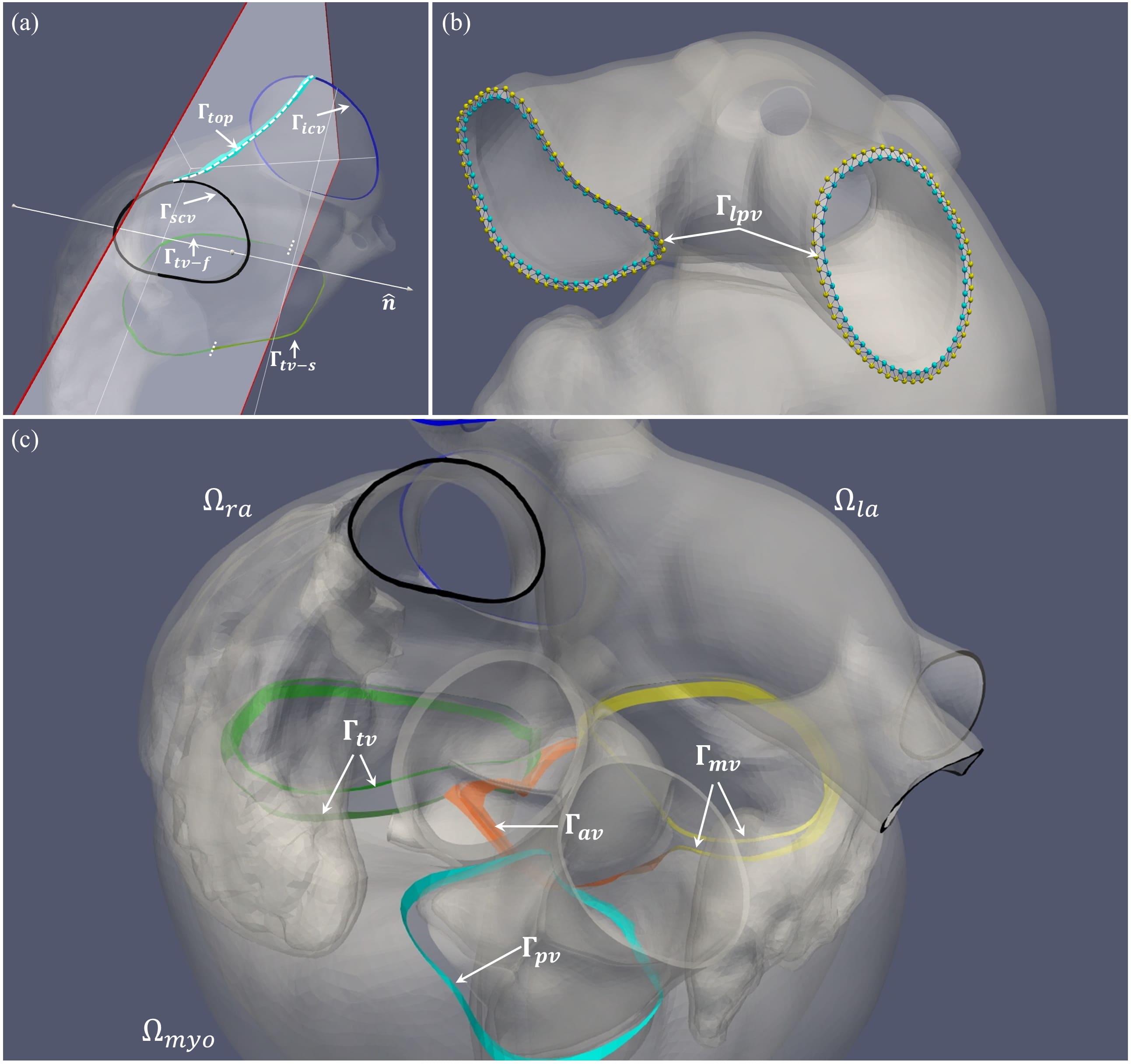}
	\caption{\textcolor{black}{Tagging procedure performed to impose the boundary conditions for the ventricular and atrial LDRBMs (presented in Sections \ref{sec:ventr} and \ref{sec:atrial}). (a) Tagging procedure for the top band $\Gamma_{top}$ and the tricuspid valve $\Gamma_{tv}$ (such that $\Gamma_{tv}=\Gamma_{tv-s} \cup \Gamma_{tv-f}$) in the right atrium: $\hat{n}$ is the normal of the plane passing trough $\Gamma_{top}$ while $\Gamma_{scv}$ and $\Gamma_{icv}$ are the rings of the caval veins; (b)~tagging procedure for the left pulmonary vein rings $\Gamma_{lpv}$ in the left atrium: yellow and blue points lay on the epicardial and endocardial border zone, respectively; (c) result of the tagging procedure for the four valve rings ($\Gamma_{mv}$ mitral valve, $\Gamma_{av}$ aortic valve, $\Gamma_{tv}$ tricuspid valve and $\Gamma_{pv}$ pulmonary valve) in a complete biventricular geometry. For further details refer to~\cite{fedele2019processing}.}}
	\label{fig:Tag}
\end{figure}


\subsection*{\textcolor{black}{Tagging procedure for the atrial and ventricular LDRBMs}}\label{sec:tag}
\textcolor{black}{For the atria the first tagging step consists in extracting the endocardium and the epicardium from the untagged surface model. Then, tags of the pulmonary and caval veins, coronary sinus, tricuspid and mitral valves rings are performed by connecting the points laying on the border zone of the endocardium to the corresponding epicardial points, see Figure \ref{fig:Tag}(b) (for further details about the connection procedure refer to \cite{fedele2019processing}). Furthermore, the tag $\Gamma_{top}$ in RA is carried out by manually producing a straightforward band connecting the top upper elements of SCV and ICV, see Figure \ref{fig:Tag}(a). Finally, the tricuspid valve ring $\Gamma_{tv}$ in RA is subdivided in one part facing the atrial septum $\Gamma_{tv-s}$ and another one related to the free wall $\Gamma_{tv-f}$, such that $\Gamma_{tv}=\Gamma_{tv-s} \cup \Gamma_{tv-f}$: this subdivision is produced by clipping the tricuspid valve ring $\Gamma_{tv}$ with a plane passing through $\Gamma_{top}$ band, see Figure \ref{fig:Tag}(a).} 

\textcolor{black}{For the ventricles, the first tagging step consists in extracting the epicardium and the right and left endocardia from the untagged surface model. Moreover, for R-RBM, the right endocardium $\Gamma_{rv}$ is subdivided into the right septum $\Gamma_{rs}$ and the remaining part $\Gamma_{rv-s}$ such that $\Gamma_{rv}=\Gamma_{rs} \cup \Gamma_{rv-s}$: this tagging subdivision is achieved by selecting a threshold in the distance between right and left endocardia, see step~1 in Figure~\ref{fig:Rossi-LDRBM}. Furthermore, concerning a based biventricular geometry, the final tagging step consists in producing an upper basal plane between the ventricular endocardium and epicardium, see step 1 in Figure~\ref{fig:Rossi-LDRBM}. Regarding a complete biventricular model, tags of  the four valve rings ($\Gamma_{mv}$ mitral valve, $\Gamma_{av}$ aortic valve, $\Gamma_{tv}$ tricuspid valve and $\Gamma_{pv}$ pulmonary valve) are defined by selecting a threshold in the distance from the corresponding atrial rings for $\Gamma_{mv}$ and $\Gamma_{tv}$ and from the aortic and pulmonary OT roots for $\Gamma_{av}$ and $\Gamma_{pv}$, respectively (see Figure \ref{fig:Tag}(c)).}    

\section*{Acknowledgements}
This project has received funding from the European Research Council (ERC) under the European Union’s Horizon 2020 research and innovation programme (grant agreement No 740132, iHEART - An Integrated Heart Model for the simulation of the cardiac function, P.I. Prof. A. Quarteroni). CV has been partially supported by the H2020-MSCA-ITN-2017, EU project 765374 "ROMSOC - Reduced Order Modelling, Simulation and Optimization of Coupled systems". LD, CV and AQ have been partially supported by the Italian research project MIUR PRIN17 2017AXL54F. "Modeling the heart across the scales: from cardiac cells to the whole organ". We are particularly grateful to the two anonymous referees for their helpful comments and suggestions.
\begin{center}
	\raisebox{-.5\height}{\includegraphics[width=.15\textwidth]{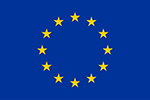}}
	\hspace{2cm} 
	\raisebox{-.5\height}{\includegraphics[width=.13\textwidth]{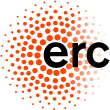}}
\end{center}


\bibliographystyle{unsrt} 
\bibliography{mylatexbi}




\end{document}